\documentclass[a4paper,11pt]{article}
\usepackage{amsfonts}
\usepackage{amssymb}
\usepackage{amsmath,amsthm}

\newtheorem{theoreme}{Th\'eor\`eme}[section]

\newtheorem{prop}[theoreme]{Proposition}
\newtheorem{cor}[theoreme]{Corollaire}
\newtheorem{definition}[theoreme]{D\'efinition}

\def \dem {\ \emph{D\'emonstration.} }
\def \cqfd {\ { \hspace{\stretch{9} } $\square$}}
\def \Rq {\ {\bf Remarque.} }

\def \sm {\setminus }
\newcommand{\be}{\begin{enumerate}}  \newcommand{\ee}{\end{enumerate}}
\newcommand{\bi}{\begin{itemize}}  \newcommand{\ei}{\end{itemize}}
\newcommand{\bd}{\begin{description}}  \newcommand{\ed}{\end{description}}
\newcommand{\ra}{\rightarrow}
\newcommand{\lan}{\langle}   \newcommand{\ran}{\rangle}

\def \R {\mathbb{R}}    \def \T {\mathbb{T}}

\renewcommand{\phi}{\varphi}
\renewcommand{\epsilon}{\varepsilon}
\title{Sur la $L^2-$cohomologie des vari\'et\'es \`a courbure n\'egative}
\author{Nader Yeganefar}
\date{}
\begin{document}
\maketitle

\begin{abstract}
We give a topological interpretation of the space of $L^2$-harmonic forms of finite-volume manifolds with sufficiently pinched negative curvature. We give examples showing that this interpretation fails if the curvature is not sufficiently pinched and that our result is sharp with respect to the pinching constants. The method consists first in comparing $L^2-$cohomology with weighted $L^2-$cohomology thanks to previous works done by T. Ohsawa, and then in identifying these weighted spaces.
\vskip0.2cm
\centerline{\bf R\'esum\'e}
\vskip0.2cm
Nous donnons une interpr\'etation topologique des espaces de formes harmoniques $L^2$ des vari\'et\'es de volume fini, \`a courbure n\'egative suffisamment pinc\'ee. Nous donnons des exemples montrant que cette interpr\'etation n'est plus valable si la courbure n'est pas suffisamment pinc\'ee et que le r\'esultat est optimal. La m\'ethode utilis\'ee consiste \`a comparer $L^2-$cohomologie et $L^2-$cohomologie \`a poids gr\^ace \`a des travaux de T. Ohsawa, puis \`a identifier la $L^2-$cohomologie \`a poids.
\end{abstract}


\section*{Introduction}

Soit $(M^n,g)$ une vari\'et\'e riemannienne compl\`ete. On note $\mathcal{H}^k(M)$ l'espace des $k-$formes harmoniques $L^2$ de $M$, c'est-\`a-dire celles qui sont de carr\'e int\'egrable, ferm\'ees et coferm\'ees. Lorsque $M$ est compacte, on sait, gr\^ace au th\'eor\`eme de Hodge-de Rham, que cet espace est de dimension finie, et isomorphe au $k-$\`eme espace de cohomologie r\'eelle de $M$. Quand $M$ est non compacte, ce que nous supposerons toujours dans la suite, il est naturel de se demander ce qui subsiste de ce r\'esultat : est-ce que l'espace $\mathcal{H}^k(M)$ est de dimension finie, et si oui peut-on en donner une interpr\'etation topologique ? Par exemple, d'apr\`es J. Lott \cite{L3}, si $M$ est une vari\'et\'e compl\`ete de volume fini, \`a courbure sectionnelle $K$ n\'egative et pinc\'ee (i.e. il existe des r\'eels $0<a\leq b$, avec $-b^2\leq K\leq -a^2$), alors tous ses espaces de formes harmoniques sont de dimension finie. De plus, J. Lott a pos\'e la question suivante (MSRI, printemps 2001) :\\
{\bf Question. }\textit{Soit $(M^n,g)$ une vari\'et\'e riemannienne compl\`ete de dimension $n$, de volume fini, \`a courbure sectionnelle n\'egative et pinc\'ee. Est-il vrai qu'on a les isomorphismes :
$$\mathcal{H}^k(M)\simeq \left\{ \begin{array}{lll}
     H^k(M), & \textrm{si $k<(n-1)/2$,}\\
{\rm Im} (H^{k}_c(M)\ra H^{k}(M)), & \textrm{si $k=n/2,(n\pm 1)/2$,}\\
     H^k_c(M), & \textrm{si $k> (n+1)/2$ ?}
     \end{array} \right.  $$}
La r\'eponse \`a cette question est "oui" dans le cas des vari\'et\'es hyperboliques de volume fini \cite{Z}, \cite{M-P} et des quotients de volume fini de l'espace hyperbolique complexe \cite{Z}. Notre r\'esultat principal est le suivant :
\begin{theoreme}\label{cohomologie intro}
Soit $(M^{n},g)$, une vari\'et\'e riemannienne compl\`ete de dimension $n$,
 de volume fini, \`a courbure sectionnelle $K$ n\'egative pinc\'ee : il existe des constantes $a$ et $b$, telles que l'on ait $-b^{2}\leq K\leq -a^{2} <0$. On suppose que $na-(n-2)b>0,$ alors on a les isomorphismes entre espaces vectoriels de dimension finie :
$$\mathcal{H}^k(M)\simeq \left\{ \begin{array}{lll}
     H^k(M), & \textrm{si $k<(n-1)/2$,}\\
{\rm Im} (H^{n/2}_c(M)\ra H^{n/2}(M)), & \textrm{si $k=n/2$,}\\
     H^k_c(M), & \textrm{si $k> (n+1)/2.$}
     \end{array} \right. $$
Si de plus la courbure est constante ($a=b$), et la dimension est impaire, alors on a $\mathcal{H}^{(n\pm 1)/2}(M)\simeq {\rm Im} (H^{(n\pm 1)/2}_c(M)\ra H^{(n\pm 1)/2}(M))$. 
\end{theoreme}
Ceci signifie que la r\'eponse \`a la question est affirmative si la courbure est suffisamment pinc\'ee. Nous \'etudierons des exemples qui montrent que la r\'eponse est n\'egative si l'on n'impose pas de conditions sur le pincement, et que notre r\'esultat est optimal.

Avant d'indiquer le sch\'ema de preuve de ce th\'eor\`eme, nous avons besoin de pr\'esenter quelques r\'esultats. Tout d'abord, il est connu que les espaces de formes harmoniques admettent une interpr\'etation en termes de $L^2-$cohomologie, avec laquelle il est souvent pr\'ef\'erable de travailler. Des rappels sur la $L^2-$cohomologie seront faits dans la premi\`ere partie. D'une part, un fait important, d\^u \`a J. Lott \cite{L2}, est que la finitude de la dimension des espaces de $L^2-$cohomologie ne d\'epend que de la g\'eom\'etrie \`a l'infini : les espaces de $L^2-$cohomologie de deux vari\'et\'es isom\'etriques en dehors d'un compact sont simultan\'ement de dimension finie ou infinie. On peut donc esp\'erer trouver des liens entre la topologie de la vari\'et\'e, la g\'eom\'etrie \`a l'infini, et la $L^2-$cohomologie. D'autre part, un cas simple o\`u on sait que l'espace des $k-$formes harmoniques est de dimension finie est celui o\`u z\'ero n'est pas dans le spectre essentiel du laplacien de Hodge-de Rham $\Delta _k$ agissant sur les $k-$formes diff\'erentielles de carr\'e int\'egrable. Or, on sait, gr\^ace aux travaux de Glazman \cite{G}, Donnelly \cite{Do}, Anghel \cite{Ang} ou encore B\"ar \cite{Ba}, que le spectre essentiel essentiel ne d\'epend que de la g\'eom\'etrie \`a l'infini. En effet, pla\c cons-nous dans un cadre g\'en\'eral o\`u $(X,g)$ est une vari\'et\'e, \'eventuellement \`a bord compact, qui est m\'etriquement compl\`ete, et soit $\Delta _k$ le laplacien d\'efini avec conditions absolues (ou relatives) au bord. Alors $0$ n'est pas dans le spectre essentiel de $\Delta _k$ si et seulement si on a une in\'egalit\'e de Poincar\'e \`a l'infini : il existe un compact $K$, et il existe une constante $C>0$ tels que pour toute $k-$forme $\alpha$ lisse \`a support compact dans l'ext\'erieur de $K$, l'on ait
$$C\Vert \alpha \Vert _{L^2}\leq  \Vert \Delta _k\alpha \Vert _{L^2}.$$
Ces consid\'erations nous am\`enent \`a montrer le r\'esultat suivant, qui, sous une hypoth\`ese portant sur le spectre essentiel du laplacien, prouve l'exactitude d'une suite qui relie la cohomologie \`a support compact $H^*_c(D)$ d'un ouvert $D$ born\'e \`a bord r\'egulier de $M$, \`a la $L^2-$cohomologie $H^*_2(M)$ de $M$, et \`a la $L^2-$cohomologie absolue $H^*_2(M\sm D)$ du compl\'ementaire de $D$ :
\begin{theoreme}
Soient $(M,g)$ une vari\'et\'e riemannienne compl\`ete, et $D$ un ouvert born\'e \`a bord r\'egulier de $M$. On suppose que pour un certain entier $k$, $0$ n'est pas dans le spectre essentiel de $\Delta _k.$ Alors nous avons la suite exacte
\begin{eqnarray*} 
H^{k-1}_2(M\setminus D)\buildrel b \over \rightarrow H^k_c(D)\buildrel e \over \rightarrow H^k_2(M)\buildrel r \over \rightarrow H^k_2(M\sm D) \\ \buildrel b \over \rightarrow H^{k+1}_c(D) \buildrel e \over \rightarrow H^{k+1}_2(M) \buildrel r \over \rightarrow H^{k+1}_2(M\sm D) \buildrel b \over \rightarrow H^{k+2}_c(D).
\end{eqnarray*}
\end{theoreme}
\Rq On \'enoncera en fait un r\'esultat un peu plus g\'en\'eral dans la premi\`ere partie. Observons aussi que cette suite exacte est bien connue. En effet, si $0$ n'est pas dans le spectre essentiel de $\Delta _k$, alors la $L^2-$cohomologie et la $L^2-$cohomologie non r\'eduite sont les m\^emes en degr\'e $k$, et on sait dans ce cas, d'apr\`es les travaux de Cheeger, que la suite est exacte au rang $k$. De plus, si l'op\'erateur de Gauss-Bonnet $d+\delta$ est ``non-parabolique \`a l'infini'', et v\'erifie une hypoth\`ese suppl\'ementaire, alors la suite exacte est vraie pour tous les degr\'es (voir \cite{C1} et \cite{C4} pour plus de d\'etails). Ce qui est int\'eressant ici, c'est que l'hypoth\`ese sur le spectre essentiel en degr\'e $k$ entra\^\i ne aussi l'exactitude de la suite au rang $k+1$. Ce r\'esultat est probablement bien connu, mais n'ayant pas trouv\'e de r\'ef\'erence pr\'ecise, nous avons pr\'ef\'er\'e l'inclure ici.

On a aussi une notion d'espaces de $L^2-$cohomologie \`a poids, et ce qui pr\'ec\`ede reste valable pour ces espaces. Pour prouver le th\'eor\`eme \ref{cohomologie intro}, notre premi\`ere t\^ache sera de montrer :
\begin{prop}\label{poidsfini intro}
Soit $(M^n,g)$ une vari\'et\'e compl\`ete de volume fini, \`a courbure sectionnelle $K$ n\'egative et pinc\'ee : il existe deux constantes strictement positives $a$ et $b$ telles que $-b^2\leq K\leq -a^2<0$. Soit $U : M \ra \R$ une fonction qui sur chaque bout de $M$ est de la forme $U=Cr$, o\`u $r$ est une fonction de Busemann associ\'ee au bout, et $C$ une constante v\'erifiant $C>2a+(n-1)b$. Alors la $L^2-$cohomologie \`a poids de $M$, associ\'ee au poids $e^U$, est isomorphe \`a la cohomologie \`a support compact de $M$ : $H^*_{2,U}(M)\simeq H^*_c(M).$
\end{prop}
Indiquons rapidement comment on prouve cette proposition. Le raisonnement est d'ailleurs plus ou moins implicite dans \cite{Do-F}. D'abord, on montre que $0$ n'est pas le spectre essentiel du laplacien associ\'e au poids. Ceci repose sur la description pr\'ecise de la g\'eom\'etrie \`a l'infini des vari\'et\'es consid\'er\'ees \cite{E}, ainsi que sur une formule d'int\'egration par parties \cite{Do-X}, qui permettent d'avoir une in\'egalit\'e de Poincar\'e \`a l'infini. On a donc la suite exacte pr\'ec\'edente pour tout ouvert born\'e $D$ (la suite exacte est aussi valable pour la $L^2$-cohomologie \`a poids). Ensuite, la m\^eme formule d'int\'egration par parties permet de montrer l'annulation de la $L^2-$cohomologie \`a poids \`a l'infini (i.e. celle de $M\sm D$) si $D$ est un ouvert bien choisi. La suite exacte permet alors de conclure.

Pour montrer le th\'eor\`eme \ref{cohomologie intro}, on doit alors comparer la $L^2-$cohomologie de $M$ \`a des espaces de $L^2-$cohomologie \`a poids, avec des poids bien choisis (voir la proposition \ref{amelioration}). Ceci se fait gr\^ace \`a des travaux d'Ohsawa \cite{O}, que l'on modifie l\'eg\`erement, et on obtient ainsi la conclusion du th\'eor\`eme pour les degr\'es $k$ loin du degr\'e moiti\'e $n/2$.  Pour les degr\'es proches du degr\'e moiti\'e, on aura besoin d'autres arguments, et notamment du th\'eor\`eme d'indice de Ballmann et Br\"uning \cite{B-B} lorsque $n$ est pair. Ceci sera pr\'ecis\'e dans la partie 5.

Le plan de cet article est le suivant : dans la premi\`ere partie, on commence par faire des rappels sur la $L^2-$cohomologie. Ensuite, on \'enonce et on d\'emontre le th\'eor\`eme concernant la suite exacte. Dans la deuxi\`eme partie, on pr\'esente la formule d'int\'egration par parties de Donnelly et Xavier (voir \cite{Do-X}, \cite{Do-F}, \cite{E-F} et \cite{B-B}), en remarquant qu'elle s'adapte aussi au cas o\`u l'on a un poids. Dans la troisi\`eme, on montre un r\'esultat concernant la $L^2-$cohomologie \`a poids. Ahmed et Stroock \cite{A-S} ont trait\'e de tels exemples, mais notre m\'ethode permet d'obtenir des r\'esultats plus g\'en\'eraux :
\begin{theoreme}
Soit $(M^n,g)$ une vari\'et\'e riemannienne compl\`ete de dimension $n$. On suppose qu'il existe une fonction $U : M \ra \R$ qui v\'erifie les hypoth\`eses suivantes :
\be
 \item $U$ exhauste $M$,
 \item il existe une constante $\varepsilon >0$ telle qu'en dehors d'un compact, on ait $$\varepsilon ^2 \leq \vert \nabla U \vert ^2 ,$$
 \item il existe deux constantes $C_1\leq 0$ et $C_2\in \R$ telles qu'en dehors d'un compact, on ait $$C_1 \vert \nabla U \vert ^2 \leq HU \leq C_2\vert \nabla U \vert ^2.$$  
\ee
Si $1+nC_1-2nC_2>0$, alors $0$ n'est pas dans le spectre essentiel de l'op\'erateur $d+\delta _U$. De plus, les espaces de $L^2$-cohomologie \`a poids de $M$, associ\'es au poids $e^U$, sont de dimension finie, et isomorphes aux espaces de cohomologie \`a support compact de $M$ : $H^*_{2,U}(M)\simeq H^*_c(M).$
\end{theoreme}
Dans la quatri\`eme partie, on \'enonce les r\'esultats d'Ohsawa, en montrant comment ses arguments s'adaptent \`a notre cas. Dans la cinqui\`eme partie, on montre le r\'esultat principal concernant la $L^2-$cohomologie des vari\'et\'es de volume fini, \`a courbure n\'egative suffisamment pinc\'ee. On poursuit par un contre-exemple montrant que le r\'esultat est faux si la courbure n'est pas assez pinc\'ee. La partie 6 contient d'autres exemples qui peuvent \^etre trait\'es par la m\'ethode qui a permis de prouver la proposition \ref{poidsfini intro}. Ainsi, on consid\`ere d'abord certains produits tordus, puis on retrouve de mani\`ere plus simple le r\'esultat de Mazzeo \cite{M} concernant les formes harmoniques des vari\'et\'es conform\'ement compactes. 

\textit{Remerciements.} Ce travail a \'et\'e effectu\'e alors que je suis en th\`ese sous la direction de G. Carron. Je tiens donc \`a remercier ce dernier de m'avoir guid\'e et conseill\'e, m'exposant ses id\'ees avec enthousiasme et toujours pr\^et \`a l'\'ecoute. Je voudrais aussi remercier tous les membres du d\'epartement de math\'ematiques de l'Universit\'e de Nantes, et en particulier les autres doctorants.
\section{Une suite exacte en $L^2-$cohomologie}

 \subsection{Rappels sur la $L^2-$cohomologie}
Nous rappelons d'abord des faits bien connus sur la $L^2-$cohomologie (voir entre autres \cite{L1}, \cite{L2}, ou \cite{C1}, \cite{C4}). Soit $(M,g)$ une vari\'et\'e riemannienne compl\`ete. On note $C^\infty _0(\Lambda ^kT^*M)$ (respectivement $L^2(\Lambda ^kT^*M)$, etc...) l'ensemble des $k-$formes lisses \`a support compact (respectivement de carr\'e int\'egrable, etc...) dans $M$. Le $k-$\`eme espace de $L^2-$cohomologie (r\'eduite) de $M$ est d\'efini par$$H^k_2(M)=\{\alpha \in L^2(\Lambda ^kT^*M)/\, d\alpha =0\}/\overline{dC^\infty _0(\Lambda ^{k-1}T^*M)}^{L^2} .$$ 
Un autre espace tr\`es proche souvent consid\'er\'e est l'espace de $L^2$-cohomologie non r\'eduite, qui est le quotient de $\{\alpha \in L^2(\Lambda ^kT^*M)/\, d\alpha =0\}$ par $\{d\alpha/\,\alpha \in L^2(\Lambda ^{k-1}T^*M),\, d\alpha \in L^2\}$, sans prendre d'adh\'erence. En g\'en\'eral, $L^2-$cohomologie r\'eduite et non r\'eduite ne sont pas \'egales. Il y a n\'eanmoins \'egalit\'e en degr\'e $k$ lorsque $0$ n'est pas dans le spectre essentiel du laplacien $\Delta _k$, et la proposition \ref{exact} montre que dans ce cas il y a aussi \'egalit\'e en degr\'e $k+1$. Dans la suite, ``$L^2$-cohomologie'' voudra dire ``$L^2$-cohomologie r\'eduite''.

Il y a une interpr\'etation de la $L^2-$cohomologie en termes de formes harmoniques. En effet, notons $\mathcal{H}^k(M)$ l'espace des $k-$formes harmoniques $L^2$ de $M$ : $$\mathcal{H}^k(M)=\{\alpha \in L^2(\Lambda ^kT^*M)/\, d\alpha =\delta \alpha =0\},$$ o\`u $\delta$ est l'op\'erateur d\'efini initialement sur les formes lisses \`a support compact comme l'adjoint de $d$. Comme $M$ est compl\`ete, $\mathcal{H}^*(M)$ est aussi le noyau $L^2$ du laplacien. Un fait important est la d\'ecomposition de Hodge-de Rham-Kodaira \cite[th\'eor\`eme 24]{dR} :
$$L^2(\Lambda ^kT^*M)= \mathcal{H} ^k(M) \oplus
  \overline{dC^\infty _0(\Lambda ^{k-1} T^*M)} \oplus
  \overline{\delta C^\infty _0(\Lambda ^{k+1} T^*M)},$$
et de plus,
  $$\{\alpha \in L^2(\Lambda ^kT^*M)/\, d\alpha =0\}=\mathcal{H} ^{k}(M) \oplus
  \overline{dC^{\infty}_{0}(\Lambda ^{k-1} T^*M)}.$$
On en d\'eduit que $$H^k_2(M)\simeq \mathcal{H}^k(M).$$
Maintenant, si $(M,g)$ est une vari\'et\'e avec un bord compact, et m\'etriquement compl\`ete, on peut aussi d\'efinir des espaces de $L^2$-cohomologie absolue et relative. On note par $C^\infty _b(\Lambda ^kT^*M)$ l'espace des $k$-formes lisses \`a support born\'e sur $M$, le support pouvant rencontrer le bord (contrairement \`a ce qui se passe pour les \'el\'ements de $C^\infty _0$). L'espace de $L^2$-cohomologie absolue est alors d\'efini par :
$$H^k_2(M)=(\delta C^\infty _0(\Lambda ^{k+1}T^*M))^\perp/\overline{dC^\infty _b(\Lambda ^{k-1}T^*M)} .$$
Cet espace est isomorphe \`a un espace de formes harmoniques v\'erifiant la condition absolue sur le bord $\partial M$ : si $i_\nu \alpha$ d\'esigne le produit int\'erieur de la forme $\alpha$ par la normale $\nu$ sur $\partial M$, alors
$$H^k_2(M)\simeq \mathcal{H}^k_A(M):=\{\alpha \in L^2(\Lambda ^kT^*M)/\, d\alpha =\delta \alpha =0,\, i_\nu\alpha =0\} .$$
La $L^2$-cohomologie relative, quant \`a elle, est d\'efinie par :
$$H^k_2(M, \partial M)=(\delta C^\infty _b(\Lambda ^{k+1}T^*M))^\perp/\overline{dC^\infty _0(\Lambda ^{k-1}T^*M)} .$$
Elle est isomorphe \`a un espace de formes harmoniques v\'erifiant la condition relative au bord : si $i$ est l'inclusion $\partial M\hookrightarrow M$, alors
$$H^k_2(M,\partial M)\simeq \mathcal{H}^k_R(M):=\{\alpha \in L^2(\Lambda ^kT^*M)/\, d\alpha =\delta \alpha =0,\, i^*\alpha =0\} .$$
Nous terminons cette partie par la pr\'esentation de la $L^2-$cohomologie \`a poids. Si $U :M\to \R$ est une fonction lisse d\'efinie sur $M$, on pose $\rm{d}v_U=e^U\rm{d}v_g$, o\`u $\rm{d}v_g$ est la mesure riemannienne associ\'ee \`a $g$. On peut alors consid\'erer l'adjoint $\delta _U$ de $d$ par rapport \`a cette mesure $\rm{d}v_U$ :
$$\delta _U=e^{-U}\delta e^U .$$
Tout ce qui a \'et\'e dit pr\'ec\'edemment reste alors valable en rempla\c cant $\delta$ par $\delta _U$. En particulier, on peut d\'efinir les espaces de $L^2$-cohomologie \`a poids $H^k_{2,U}(M)$ par rapport au poids $U$ (ou plus pr\'ecis\'ement par rapport au poids $e^U$). De plus, si on note $L^2_U$ l'espace des formes de carr\'e int\'egrable par rapport \`a la mesure $\rm{d}v_U$, alors on a l'identification $$H^k_{2,U}(M)\simeq \mathcal{H}^k_U(M):=\{\alpha \in L^2_U(\Lambda ^kT^*M)/\, d\alpha =\delta _U \alpha =0\}.$$

 \subsection{$L^2-$cohomologie et spectre du laplacien} 
{\bf Remarque pr\'eliminaire.} Dans tout ce qui suit, nous ne consid\'erons, pour simplifier, que la $L^2$-cohomologie sans poids. Cependant, tous les r\'esultats restent aussi valables avec un poids.

Dans cette partie, nous pr\'esentons la suite exacte dont il est question dans l'introduction.
\begin{theoreme}\label{suite exacte}
Soient $(M,g)$ une vari\'et\'e riemannienne compl\`ete, et $D$ un ouvert \`a bord compact et r\'egulier de $M$. On suppose que pour un certain entier $k$, $0$ n'est pas dans le spectre essentiel de $\Delta _k.$ Alors nous avons la suite exacte
\begin{eqnarray*}
H^{k-1}_2(M)\buildrel r \over \ra H^{k-1}_2(M\sm D)\buildrel b \over \ra H^k_2(D,\partial D)\buildrel e \over \ra H^k_2(M) \\ \buildrel r \over \ra H^k_2(M\sm D) \buildrel b \over \ra H^{k+1}_2(D,\partial D) \buildrel e \over \ra H^{k+1}_2(M) \buildrel r \over \ra H^{k+1}_2(M\sm D).
\end{eqnarray*}
\end{theoreme}
Avant de d\'emontrer ce th\'eor\`eme, nous avons besoin de quelques r\'esultats pr\'eliminaires, mais avant tout, pr\'ecisons les applications $r$, $e$ et $b$ qui interviennent dans la suite exacte.
\bi
 \item $r:H^{*}_2(M)\ra H^{*}_2(M\sm D)$ est l'application de restriction induite en $L^2-$cohomologie par l'inclusion $M\sm D\hookrightarrow M$.
 
 \item $e$ est l'extension par z\'ero : si $[ \alpha ] \in H^{*}_2(D,\partial D)$, de repr\'esentant $\alpha$, on d\'efinit $\tilde {\alpha}$ sur $M$ par $\tilde {\alpha}=\alpha $ sur $D$, et $\tilde {\alpha} =0$ sur $M\sm D$. On v\'erifie que $\tilde {\alpha}$ est faiblement ferm\'ee, et on pose $e[\alpha ]=[\tilde {\alpha} ].$ Ceci est ind\'ependant du repr\'esentant choisi.

 \item $b$ se d\'efinit  comme l'homomorphisme cobord ordinaire en cohomologie de de Rham : si $[\alpha ]$ est une classe dans $H^k_2(M\sm D)$, on peut choisir un repr\'esentant $\alpha$ qui soit ferm\'e et lisse. Il existe alors une $k-$forme $\overline{\alpha}$ lisse sur $M$, qui co\"incide avec $\alpha$ sur $M\sm D$, et qui v\'erifie $d\overline{\alpha} =0$ sur un voisinage de $M\sm D$. On impose de plus que $\overline{\alpha}$ et $d\overline{\alpha}$ soient de carr\'e int\'egrable (c'est toujours possible, car le bord $\partial D$ \'etant compact, on peut choisir $\overline{\alpha}$ telle que $\overline{\alpha} \vert _D$ soit \`a support born\'e). On v\'erifie que la classe de $d\overline{\alpha}$ dans $H^{k+1}_2(D, \partial D)$ est ind\'ependante du repr\'esentant lisse $\alpha$ choisi, ainsi que du prolongement $\overline{\alpha}$,  et on pose $b[\alpha]=[d\overline{\alpha}].$   
\ei
Pour la d\'emonstration du th\'eor\`eme \ref{suite exacte}, nous aurons besoin du r\'esultat classique suivant (version facile du th\'eor\`eme 1.5 de \cite{C1}) :

\begin{prop}\label{Green}
Soit $(M,g)$ une vari\'et\'e compl\`ete. On suppose que pour un entier $k$ donn\'e, $0$ n'est pas dans le spectre essentiel de $\Delta _k.$ Alors il existe un op\'erateur de Green continu
$$ G=G_{-}\oplus G_{+}\, : \, L^2(\Lambda ^{k}T^{*}M)\ra  L^2(\Lambda ^{k-1}T^{*}M)\oplus  L^2(\Lambda ^{k+1}T^{*}M), $$
tel que pour toute forme $\alpha \in  L^2(\Lambda ^{k}T^{*}M)$, on ait
$$\alpha = P\alpha +dG_{-}\alpha +\delta G_{+}\alpha,$$
o\`u $P$ d\'esigne la projection orthogonale sur le noyau du laplacien. De plus, si $\alpha$ est lisse, alors $G\alpha$ est aussi lisse. Enfin, si $\alpha$ est nulle dans $H^k_2(M)$, alors $\alpha =dG_{-}\alpha$.
\end{prop}
\dem Comme $0$ n'est pas dans le spectre essentiel de $\Delta _k,$ il existe un op\'erateur de Green (born\'e) 
$$A\, : \, L^2(\Lambda ^{k}T^{*}M)\ra \rm{Dom}(\Delta ^{k})\subset L^2(\Lambda ^{k}T^{*}M),$$ tel que $\Delta A=\rm{Id}-P$ sur $L^2(\Lambda ^{k}T^{*}M)$. L'op\'erateur $G$ se d\'efinit par
 $$G\alpha =\delta A\alpha +dA\alpha.$$
$G$ est born\'e sur $L^2$, car $A$ envoie contin\^ument $L^2$ dans le domaine de $\Delta$, et on a $(d+\delta )^{2}=\Delta$ au sens op\'erateur. On voit alors facilement que $G$ v\'erifie les assertions de l'\'enonc\'e.\cqfd \\
\Rq Comme on l'a rappel\'e dans l'introduction, l'hypoth\`ese sur le spectre essentiel ne d\'epend que de la g\'eom\'etrie \`a l'infini. Ainsi, si on revient aux hypoth\`eses de la proposition, le laplacien d\'efini sur $D$ ou $M\sm D$ avec conditions absolues ou relatives au bord  v\'erifie aussi l'assertion sur le spectre essentiel, donc la m\^eme preuve s'adapte \`a ce cas, et on obtient l'existence d'un op\'erateur de Green.

La proposition suivante est essentielle pour montrer le d\'ecalage dans la suite exacte :
\begin{prop}\label{exact}
Soit $(X,g)$ une vari\'et\'e riemannienne avec \'eventuellement un bord compact, m\'etriquement compl\`ete. On suppose que $0$ n'est pas dans le spectre essentiel du laplacien $\Delta _k$, d\'efini avec conditions absolues ou relatives au bord. Soit $\alpha$ une forme diff\'erentielle dans $L^2\cap C^\infty$, de degr\'e $k$ ou $k+1$. Si $\alpha$ est nulle dans $H^*_2(X)$ (respectivement dans $H^*_2(X,\partial X)$), alors il existe une forme $\beta$ lisse dans le domaine $\rm{Dom}((d+\delta )^A)$ de l'op\'erateur $d+\delta$ d\'efini avec condition absolue sur le bord (respectivement dans le domaine $\rm{Dom}((d+\delta )^R)$ de l'op\'erateur $d+\delta$ d\'efini avec condition relative sur le bord), coferm\'ee, et v\'erifiant $\alpha =d\beta$.
\end{prop}
\dem Le cas o\`u $\alpha$ est de degr\'e $k$ d\'ecoule de la proposition \ref{Green} et de la remarque qui suit sa preuve. On fait la d\'emonstration uniquement dans le cas o\`u $X$ a un bord, et o\`u $\alpha$ est de degr\'e $k+1$, nulle dans $H^*_2(X)$. Cette nullit\'e signifie qu'il existe une suite $\beta _l$ de $k-$formes lisses \`a supports born\'es telle que $\alpha =\lim d\beta _l$. On note $\Delta _k^A$ le laplacien d\'efini avec conditions absolues sur le bord, $P^A$ la projection orthogonale sur le noyau $\rm{Ker}(\Delta _k^A)$, et $(d+\delta )_k^A$ la restriction de l'op\'erateur $(d+\delta )^A$ aux $k-$formes. L'hypoth\`ese sur le spectre essentiel montre, gr\^ace \`a la proposition \ref{Green}, qu'il existe un op\'erateur de Green $G^A$ pour l'op\'erateur $(d+\delta )^A_k$, qui envoie l'espace $L^2$ dans le domaine de $(d+\delta )^A$. On peut donc \'ecrire
$$\beta _l=P^A(\beta _l)+d G^A_{-}\beta _l +\delta G^A_{+}\beta _l.$$
Or $G^A_{+}\beta _l$ v\'erifie la condition absolue au bord : $i_\nu G^A_{+}\beta _l=0$, o\`u $i_\nu$ d\'esigne le produit int\'erieur par la normale $\nu$ sur le bord $\partial X$ ; donc $\delta G^A_{+}\beta _l$ v\'erifie aussi cette condition absolue. Par cons\'equent, quitte \`a remplacer $\beta _l$ par $\delta G^A_{+}\beta _l$, on peut supposer que $\beta _l$ est dans le domaine de $(d+\delta )^A_k$, coferm\'ee, orthogonale au noyau $\rm{Ker}(\Delta ^A_k)$, et que les $d\beta _l$ (qui restent dans $C^\infty _b$) continuent de converger vers $\alpha$. De plus, il existe une constante $C>0$ telle que :
$$\forall \phi \in \rm{Dom}((d+\delta )_k^A), \, \phi \bot \rm{Ker}(\Delta ^A _k),\,  C\Vert \phi \Vert ^2 _{L^2}\leq \Vert d\phi \Vert ^2 _{L^2} +\Vert \delta \phi \Vert ^2 _{L^2}.$$
Cette in\'egalit\'e montre alors que $(\beta _l)$ est une suite de Cauchy dans $L^2$, par cons\'equent il existe une $k-$forme $\beta$ telle que $\beta _l\rightarrow \beta$ dans $L^2$. On a donc $d\beta _l\rightarrow d\beta$ au sens des distributions, mais comme $d\beta _l\rightarrow \alpha$ dans $L^2$, ceci montre que $\alpha =d\beta$ au sens des distributions. D'autre part, les $\beta _l$ sont coferm\'ees, donc on a aussi $\delta \beta =0$ au sens des distributions. Ainsi, on a $(d+\delta) (\beta)=\alpha$ au sens des distributions, et comme $\alpha$ est lisse, par r\'egularit\'e elliptique de l'op\'erateur $d+\delta$, $\beta$ est aussi lisse, et on a finalement $\alpha =d\beta.$ Enfin, comme $\beta _l\in\rm{Dom}((d+\delta )_k^A)$ converge vers $\alpha$ dans $L^2$, et que $(d+\delta )^A \beta _l$ converge vers $d\alpha$ dans $L^2$, le fait que l'op\'erateur $(d+\delta )^A$ soit ferm\'e entra\^\i ne que $\beta$ est dans le domaine $\rm{Dom}((d+\delta )_k^A)$.\cqfd

On va maintenant pouvoir passer \`a la preuve du th\'eor\`eme \ref{suite exacte} (voir aussi \cite{C1} ou \cite{C4}). Tout d'abord, par construction on a toujours
$$r\circ e =0,\,\, b\circ r=0, \,\, e\circ b=0,$$ donc
$$\rm{Im}(e) \subset \rm{Ker}(r),\,\, \rm{Im}(r) \subset \rm{Ker}(b),\,\, \rm{Im}(b) \subset \rm{Ker}(e).$$

Montrons que $\rm{Ker}(r)\subset\rm{Im}(e).$ Soit $[\alpha ]\in  H^*_2(M)$ telle que $r[\alpha ]=0.$ Si $\alpha$ est de degr\'e $k$ ou $k+1$, gr\^ace \`a la proposition \ref{exact}, il existe une forme lisse $\beta$ dans $L^2(\Lambda ^kT^*(M\sm D))$, telle que $\alpha \vert _{M\sm D}=d\beta$. Si $\overline{\beta } \in L^2(\Lambda ^kT^*M)$ est une extension lisse de $\beta$, alors on voit que $e[\alpha -d\overline{\beta }]=[\alpha ]$.

Montrons maintenant que $\rm{Ker}(e)\subset\rm{Im}(b).$ Soit $\alpha$ une forme de degr\'e $k$ ou $k+1$, de carr\'e int\'egrable sur $D$, avec $d\alpha =\delta \alpha =0$, et $i^{*}\alpha=0,$ o\`u $i$ est l'inclusion $\partial D \hookrightarrow D$. Si $\tilde{\alpha}$ d\'esigne l'extension par z\'ero de $\alpha$, on suppose que $\tilde{\alpha}$ est nulle en $L^2-$cohomologie (i.e. $[\alpha ]\in \rm{Ker}(e)$). Gr\^ace \`a la proposition \ref{exact}, on \'ecrit $\tilde{\alpha}=d\beta$, avec $\beta$ coferm\'ee, et comme $\tilde{\alpha}$ est lisse sur $D$ et $M\sm D$, $\beta$ est aussi lisse sur $D$ et sur $M\sm D$. Soit $\overline{\beta }$ une extension lisse de $\beta |_{M\sm D}$ comme dans la d\'efinition de l'application $b$. Alors on a $\alpha -d\overline{\beta} =d(\beta-\overline{\beta})$, avec $\beta-\overline{\beta}$ lisse sur $D$, de carr\'e int\'egrable, et de tir\'e en arri\`ere nul sur $\partial D$. Donc $\alpha -d\overline{\beta}$ est $L^2-$cohomologue \`a z\'ero, c'est-\`a-dire que $[\alpha ]=b[\beta |_{M\sm D}].$

Enfin, montrons que $\rm{Ker}(b)\subset \rm{Im}(r)$. Soit $[\alpha ]\in \rm{Ker}(b)$ : avec les notations pr\'ec\'edentes, ceci signifie que $d\overline{\alpha} | _D$ est nulle dans $H^*_2(D,\partial D)$. Si $\alpha$ est de degr\'e $k-1$ ou $k$, alors d'apr\`es la proposition \ref{exact}, il existe une forme $\beta$ d\'efinie sur $D$, de carr\'e int\'egrable, v\'erifiant la condition relative $i^*\beta =0$ au bord, et telle que $d\overline{\alpha} | _D=d\beta$. On note $\tilde{\beta }$ l'extension par z\'ero de $\beta$ sur $M$. La forme $\overline{\alpha } -\tilde{\beta }$ est ferm\'ee (faiblement), de carr\'e int\'egrable, et on a $r[\overline{\alpha } -\tilde{\beta }]=[\alpha ].$\cqfd

Dans le but de relier $L^2-$cohomologie et topologie, il est utile de remarquer :
\begin{cor}\label{se}
Soient $(M,g)$ une vari\'et\'e riemannienne compl\`ete, et $D$ un ouvert born\'e \`a bord r\'egulier de $M$. On suppose que pour un certain entier $k$, $0$ n'est pas dans le spectre essentiel de $\Delta _k.$ Alors nous avons les deux suites exactes
\begin{eqnarray*}
H^{k-1}_2(M\setminus D)\buildrel b \over \rightarrow H^k_c(D)\buildrel e \over \rightarrow H^k_2(M)\buildrel r \over \rightarrow H^k_2(M\sm D) \\ \buildrel b \over \rightarrow H^{k+1}_c(D) \buildrel e \over \rightarrow H^{k+1}_2(M) \buildrel r \over \rightarrow H^{k+1}_2(M\sm D) \buildrel b \over \rightarrow H^{k+2}_c(D) ,
\end{eqnarray*}
et
\begin{eqnarray*}
 H^{k-1}_2(M\sm D,\partial D)\buildrel e \over \rightarrow H^{k-1}_2(M)\buildrel r \over \rightarrow H^{k-1}(D) \buildrel b \over \rightarrow H^{k}_2(M\sm D,\partial D)\\ \buildrel e \over \rightarrow H^k_2(M) \buildrel r \over \rightarrow H^k(D) \buildrel b \over \rightarrow H^{k+1}_2(M\sm D,\partial D)\buildrel e \over \rightarrow H^{k+1}_2(M)\buildrel r \over \rightarrow H^{k+1}(D).
\end{eqnarray*}
\end{cor}
\dem La premi\`ere suite est une cons\'equence directe du th\'eor\`eme \ref{suite exacte}, si l'on observe d'une part que $H^*_2(D,\partial D)\simeq H^*_c(D),$ car $D$ est born\'e, et d'autre part que l'\'egalit\'e $\rm{Ker}(b)=\rm{Im}(r)$ est toujours vraie lorsque $D$ est born\'e (c'est une cons\'equence de la suite usuelle en cohomologie de de Rham, voir par exemple \cite{C4} ou \cite{L2}). La deuxi\`eme suite d\'ecoule aussi du th\'eor\`eme \ref{suite exacte}, si l'on \'echange $D$ en $M\sm D$ dans ce th\'eor\`eme, puis si l'on se sert du fait que l'hypoth\`ese $D$ born\'e entra\^\i ne $H^*_2(D)\simeq H^*(D)$, et $\rm{Ker}(r)=\rm{Im}(e).$\cqfd

\section{Une formule d'int\'egration par parties}
Le but de cette partie est de pr\'esenter une formule d'int\'egration par parties pour l'op\'erateur $d+\delta$. Cette formule et ses corollaires seront utiles pour l'\'etude de la $L^2-$cohomologie de certains types de vari\'et\'es. Le th\'eor\`eme principal est d'abord d\^u \`a Donnely et Xavier \cite{Do-X} ; une d\'emonstration plus synth\'etique en est aussi donn\'ee par Escobar et Freire dans \cite{E-F}. Ensuite Ballmann et Br\"uning \cite{B-B} ont \'etendu cette formule au cadre plus large des op\'erateurs de type Dirac et ont am\'elior\'e les r\'esultats de Donnelly et Xavier.

Avant d'\'enoncer le r\'esultat principal, nous faisons la remarque suivante : si $f$ est une fonction de classe $C^{2}$ sur une vari\'et\'e riemannienne $M$, alors son hessien $Hf$ se prolonge \`a une forme bilin\'eaire sym\'etrique sur $\Lambda ^{k}T^{*}M$ comme suit
\begin{equation}\label{hessien} 
Hf(\alpha,\beta)=\sum _{i,j}Hf(X_{i},X_{j})\langle i_{X_{i}}\alpha,i_{X_{j}}\beta \rangle,
\end{equation}
o\`u les $X_{i}$ forment un rep\`ere orthonormal local et $i_{X}\alpha$ d\'esigne le produit int\'erieur de la forme $\alpha$ par le champ de vecteurs $X$ (la formule est bien entendu ind\'ependante du choix d'un tel rep\`ere). Nous avons alors, suivant \cite{E-F} : 
\begin{theoreme}\label{IPP}
On consid\`ere $(M^{n},g)$ une vari\'et\'e riemannienne compl\`ete. Soient $f$ une fonction de classe $C^{2}$ sur $M$ de gradient born\'e,  $\Omega$ un ouvert (\`a bord compact assez r\'egulier) de $M$, et $\alpha$ une forme diff\'erentielle dans le domaine de l'op\'erateur $d+\delta$, lisse jusqu'au bord de $\Omega$. On note $\nu$ la normale unitaire sortante sur $\partial \Omega$, et $\rm{d}\sigma$ la mesure riemannienne induite sur le bord $\partial \Omega$. Alors 
\begin{multline}
 \int _{\Omega} \left[ 2Hf(\alpha,\alpha )+\Delta f \vert \alpha \vert^{2} \right]\rm{d}v_g = \\
  2\int_{\Omega} \left[\langle i_{\nabla f}d\alpha ,\alpha \rangle +\langle i_{\nabla f}\alpha ,\delta \alpha \rangle \right ] \rm{d}v_g+\int _{\partial \Omega} \left[ -\vert \alpha \vert ^{2} \frac{\partial f}{\partial \nu} +2 \langle i_{\nabla f}\alpha ,i_{\nu}\alpha \rangle \right] \rm{d}\sigma. \nonumber 
\end{multline}
\end{theoreme}
Nous ne refaisons pas la d\'emonstration. L'id\'ee est de montrer l'\'equation 
\begin{equation}\label{ponctuel}
Hf(\alpha,\alpha )+\lan \nabla _{\nabla f}\alpha ,\alpha\ran =\lan di_{\nabla f}\alpha ,\alpha \ran +\lan i_{\nabla f}d\alpha ,\alpha \ran ,
\end{equation} valable pour toute forme diff\'erentielle $\alpha$. Puis on int\`egre cette relation sur un ouvert par rapport \`a la mesure riemannienne $\rm{d}v_g$, et on fait des int\'egrations par parties sur les deux termes o\`u apparaissent $\lan \nabla _{\nabla f}\alpha ,\alpha\ran =1/2\nabla f.|\alpha |^2$ et $\lan di_{\nabla f}\alpha ,\alpha \ran$.

Le premier corollaire que nous en d\'eduisons se trouve dans \cite{Do-X} et \cite{B-B} sous une forme l\'eg\`erement diff\'erente, et consiste \`a consid\'erer le cas particulier des formes \`a support compact.
\begin{cor}\label{IPP1}
Soit $f$ une fonction de classe $C^{2}$ sur $M^{n}$ telle que $\vert \nabla f\vert$ soit born\'e sur $M$. On consid\`ere $\alpha$ une forme lisse \`a support compact, alors on a l'estimation
\begin{multline}
\vert \int _{M} \left[ 2Hf(\alpha,\alpha )+\Delta f \vert \alpha \vert^{2}
 \right] \vert \leq 2 \sup _{M} (\vert \nabla f \vert) \left[ \Vert d\alpha \Vert _{L^2}+\Vert \delta \alpha \Vert _{L^2}\right] \Vert \alpha \Vert _{L^2} \nonumber
\end{multline}
\end{cor}
\dem On consid\`ere un ouvert $\Omega$ qui contient le support de $\alpha$ (dans son int\'erieur), donc dans l'\'egalit\'e du th\'eor\`eme pr\'ec\'edent, il n'y a plus de termes de bord. On majore le deuxi\`eme membre de cette \'egalit\'e en utilisant le fait que $\vert i_{X}\alpha \vert\leq \vert X\vert \vert \alpha \vert$ pour tout champ de vecteurs $X$, puis en utilisant des in\'egalit\'es de Cauchy-Schwarz, pour arriver au r\'esultat annonc\'e.\cqfd

Cette id\'ee s'adapte aussi au cas o\`u l'on a un poids $U :M\to \R$. En effet, si on repart de l'\'egalit\'e \ref{ponctuel}, et si on l'int\`egre par rapport \`a la mesure $\rm{d}v_U=e^U\rm{d}v_g$, alors par int\'egration par parties, on trouve
\begin{cor}\label{IPPpoids}
Soient $f$ et $U$ deux fonctions de classe $C^{2}$ sur $(M,g)$, et $\Omega$ un ouvert \`a bord compact et r\'egulier de $M$. Si $\alpha$ est une forme diff\'erentielle lisse et \`a support born\'e dans $\Omega$ (le support pouvant rencontrer le bord), on a
\begin{multline} 
 \int _{\Omega} \left[ 2Hf(\alpha,\alpha )+\Delta f \vert \alpha \vert^{2} -\lan \nabla f,\nabla U\ran |\alpha |^2 \right]\rm{d}v_U =\\ 2\int_{\Omega} \left[\langle i_{\nabla f}d\alpha ,\alpha \rangle +\langle i_{\nabla f}\alpha ,\delta \alpha \rangle \right ]\rm{d}v_U +\int _{\partial \Omega} \left[ -\vert \alpha \vert ^{2} \frac{\partial f}{\partial \nu} +2 \langle i_{\nabla f}\alpha ,i_{\nu}\alpha \rangle \right]e^U\rm{d}\sigma. \nonumber 
\end{multline} 

\end{cor}


\section{Un cas de $L^2-$cohomologie \`a poids}

 \subsection{Cas "g\'en\'eral"}
 Soient $(M^n,g)$ une vari\'et\'e riemannienne non compacte, et $U : M \ra [0,\infty  [$ une fonction qui cro\^\i t assez vite \`a l'infini pour que le volume de $M$ par rapport \`a la mesure $\rm{d}v_{(-U)}=e^{-U}\rm{d}v_g$ soit fini, o\`u $\rm{d}v_g$ est la mesure riemannienne associ\'ee \`a $g$. Dans \cite{A-S}, Ahmed et Stroock montrent, sous certaines hypoth\`eses g\'eom\'etriques (courbure de Ricci minor\'ee, op\'erateur de courbure major\'e), et d'autres hypoth\`eses portant sur la fonction $U$, que les espaces de cohomologie $L^2$ \`a poids sont de dimension finie, et isomorphes aux espaces de cohomologie r\'eelle de $M$. Nous allons voir que sous des hypoth\`eses beaucoup moins restrictives, ce r\'esultat reste valable. En fait, nous nous int\'eressons d'abord \`a des poids $e^U$, plut\^ot que $e^{-U}$, mais la preuve de notre r\'esultat s'applique indiff\'eremment aux deux cas. 
\begin{theoreme}\label{poids}
Soit $(M^n,g)$ une vari\'et\'e riemannienne compl\`ete de dimension $n$. On suppose qu'il existe une fonction $U : M \ra \R$ qui v\'erifie les hypoth\`eses suivantes :
\be
 \item $U$ exhauste $M$,
 \item il existe une constante $\varepsilon >0$ telle qu'en dehors d'un compact, on ait $$\varepsilon ^2 \leq \vert \nabla U \vert ^2 ,$$
 \item il existe deux constantes $C_1\leq 0$ et $C_2\in \R$ telles qu'en dehors d'un compact, on ait $$C_1 \vert \nabla U \vert ^2 \leq HU \leq C_2\vert \nabla U \vert ^2.$$  
\ee
Si $\eta :=1+nC_1-2nC_2>0$, alors $0$ n'est pas dans le spectre essentiel de l'op\'erateur $d+\delta _U$. De plus, les espaces de $L^2$-cohomologie \`a poids de $M$, associ\'es au poids $e^U$, sont de dimension finie, et isomorphes aux espaces de cohomologie \`a support compact de $M$ : $H^*_{2,U}(M)\simeq H^*_c(M).$
\end{theoreme}
\Rq Si $C_2\leq 0$, alors on a le m\^eme r\'esultat en supposant que $1+nC_1>0$ (voir la preuve du th\'eor\`eme).\\
\dem On montre d'abord que sous nos hypoth\`eses, z\'ero n'est pas dans le spectre essentiel de $d+\delta _U$. Soit $c>0$ un r\'eel tel qu'en dehors de l'ouvert born\'e $D=\{ U<c\}$, les hypoth\`eses 2 et 3 du th\'eor\`eme soient v\'erifi\'ees. Si $\alpha$ est une forme de degr\'e $k$ sur $M\sm D$, lisse et \`a support born\'e, on a, d'apr\`es la formule d'int\'egration par parties du corollaire \ref{IPPpoids} (en prenant $f=U$ dans les notations de ce corollaire) :
\begin{multline}\label{IPP2} 
 \int _{M\sm D} \left[ 2HU(\alpha,\alpha )+\Delta U \vert \alpha \vert^{2} -|\nabla U |^2|\alpha |^2 \right]\rm{d}v_U =\\ 2\int_{M\sm D} \left[\langle i_{\nabla U}d\alpha ,\alpha \rangle +\langle i_{\nabla U}\alpha ,\delta \alpha \rangle \right ]\rm{d}v_U +\int _{\partial D} \left[ -\vert \alpha \vert ^{2} \frac{\partial U}{\partial \nu} +2 \langle i_{\nabla U}\alpha ,i_{\nu}\alpha \rangle \right]e^U\rm{d}\sigma , 
\end{multline} 
o\`u $\nu$ est la normale unitaire ext\'erieure sur $\partial D$, donc $\nu =-\nabla U/|\nabla U |$ ici. Or, on a l'estimation suivante :
\begin{eqnarray}
2HU(\alpha,\alpha )+\Delta U \vert \alpha \vert^{2} -|\nabla U |^2 |\alpha |^2 &\leq& (-1+2kC_2-nC_1)|\nabla U|^2 |\alpha |^2 \nonumber \\
&\leq& -\eta |\nabla U|^2 |\alpha |^2. \nonumber
\end{eqnarray}
On consid\`ere de plus que $\alpha$ est une forme diff\'erentielle v\'erifiant la condition absolue au bord $i_\nu\alpha =0$. Le terme de bord de l'\'egalit\'e \ref{IPP2} est alors positif car $-\partial U/\partial \nu=|\nabla U|$, donc compte tenu de l'in\'egalit\'e pr\'ec\'edente, on a  
\begin{eqnarray}
\eta \int _{M\sm D} |\nabla U|^2 |\alpha |^2 \rm{d}v_U &\leq &-2\int_{M\sm D} \left[\langle i_{\nabla U}d\alpha ,\alpha \rangle +\langle i_{\nabla U}\alpha ,\delta _U \alpha \rangle \right ] \rm{d}v_U \nonumber \\
 &\leq & 2\left( \int _{M\sm D} |\nabla U|^2 |\alpha |^2 \rm{d}v_U \right) ^{1/2} \nonumber \\
 & & \left[ ||d\alpha ||_{L^2_U(M\sm D)} +||\delta _U \alpha ||_{L^2_U(M\sm D)} \right] , \nonumber
\end{eqnarray}
o\`u pour passer de la premi\`ere \`a la deuxi\`eme in\'egalit\'e, on a d'abord utilis\'e l'in\'egalit\'e $|i_X\beta |\leq |X||\beta |$, valable pour tout champ de vecteurs $X$ et toute forme $\beta$, et ensuite des in\'egalit\'es de Cauchy-Schwarz. On en d\'eduit finalement que pour toute forme $\alpha$ \`a support born\'e dans $M\sm D$, v\'erifiant la condition absolue $i_\nu \alpha =0$ :
\begin{equation}\label{PI}
\quad \eta \varepsilon /2 ||\alpha ||_{L^2_U(M\sm D)} \leq ||d\alpha ||_{L^2_U(M\sm D)} +||\delta _U \alpha ||_{L^2_U(M\sm D)}
\end{equation}
Cette in\'egalit\'e est en particulier vraie pour des formes lisses \`a support compact dans $M\sm D$, donc on a une in\'egalit\'e de Poincar\'e \`a l'infini, et par cons\'equent $0$ n'est pas dans le spectre essentiel de l'op\'erateur $d+\delta _U$. On a alors une suite exacte pour la $L^2-$cohomologie \`a poids, analogue \`a la premi\`ere suite exacte du corollaire \ref{se}, c'est-\`a-dire que pour tout entier $k$, la suite 
$$H^k_c(D)\\ \buildrel e \over \rightarrow H^k_{2,U}(M) \buildrel r \over \rightarrow H^k_{2,U}(M\sm D) \buildrel b \over \rightarrow H^{k+1}_c(D)$$ est exacte.
On va montrer que la $L^2-$cohomologie absolue \`a poids de $M\sm D$ est nulle (c'est-\`a-dire $H^*_{2,U}(M\sm D)=0$). Soit donc $\alpha$ une forme de $L^2_U(M\sm D)$, $(d+\delta _U)$-harmonique, et satisfaisant la condition absolue au bord $i_\nu\alpha =0$. On consid\`ere une suite $\rho _l$ de fonctions lisses d\'efinies sur $M\sm D$, \`a supports born\'es, v\'erifiant :
$$|\rho _l|\leq 1,\quad \rho _l\ra 1 \,(l\ra\infty ),\quad |d\rho _l|\leq 1,\quad \textrm{et $d\rho _l\ra 0$ sur les compacts}.$$
On applique alors l'in\'egalit\'e \ref{PI} \`a $\rho _l\alpha$, puis on fait tendre $l$ vers l'infini, pour montrer que $\alpha$ est nulle. Ceci montre que $$H^k_{2,U}(M)\simeq H^k_c(D),$$ et ceci pour tout ouvert $D=\{ U<c\}$. En faisant tendre $c$ vers l'infini, on trouve le r\'esultat souhait\'e. \cqfd 

 \subsection{Cas des vari\'et\'es de volume fini, \`a courbure n\'egative et pinc\'ee}
Dans cette partie, $(M^{n},g)$ est une vari\'et\'e riemannienne compl\`ete de volume fini. On suppose de plus que la courbure sectionnelle $K$ de $M$ est n\'egativement pinc\'ee, c'est-\`a-dire qu'il existe des constantes $a$ et $b$, avec $$-b^{2}\leq K\leq -a^{2} <0.$$ On rappelle des r\'esultats concernant la g\'eom\'etrie de telles vari\'et\'es (voir \cite{E}, \cite{H-I}, ou encore le r\'esum\'e de \cite{L3}) : $M$ est diff\'eomorphe \`a l'int\'erieur d'une vari\'et\'e compacte \`a bord. Si $\Sigma$ est l'une des composantes (en nombre fini) du bord, il y a un bout correspondant. A tout rayon g\'eod\'esique d'un bout, on associe une fonction de Busemann $r$ (deux telles fonctions ne diff\`erent que d'une constante). Dans le cas g\'en\'eral, $r$ est seulement de classe $C^2$. Si $B$ est un bout associ\'e \`a la composante de bord $\Sigma$, alors $B$ est $C^2-$diff\'eomorphe \`a  $]0,\infty [ \times \Sigma$, les tranches $\{ t\} \times \Sigma$ \'etant les niveaux d'une fonction de Busemann $r$. Par ce diff\'eomorphisme, la m\'etrique s'\'ecrit $$g=dr^2+g_r,$$ o\`u les $g_r$ sont une famille de m\'etriques sur $\Sigma$.

La preuve du th\'eor\`eme \ref{poids} s'adapte aussi au cas des vari\'et\'es de volume fini, \`a courbure n\'egative et pinc\'ee :
\begin{cor}\label{poidsfini}
Soit $(M^n,g)$ une vari\'et\'e compl\`ete de volume fini, \`a courbure sectionnelle $K$ n\'egative et pinc\'ee : il existe deux constantes strictement positives $a$ et $b$ telles que $-b^2\leq K\leq -a^2<0$. Soit $U : M \ra \R$ une fonction qui sur chaque bout de $M$ est de la forme $U=Cr$, o\`u $r$ est une fonction de Busemann associ\'ee au bout, et $C$ une constante v\'erifiant $C>2a+(n-1)b$. Alors la $L^2-$cohomologie \`a poids de $M$, associ\'ee au poids $e^U$, est isomorphe \`a la cohomologie \`a support compact de $M$ : $H^*_{2,U}(M)\simeq H^*_c(M).$
\end{cor}
\dem On proc\`ede de mani\`ere similaire \`a la preuve pr\'ec\'edente (th\'eor\`eme \ref{poids}), en majorant le membre de gauche de la formule d'int\'egration par parties \ref{IPP2}. Pour cela, on se sert du th\'eor\`eme de comparaison des hessiens de Greene et Wu (voir \cite[p. 19]{G-W}) qui montre que si $\lambda _1\geq \ldots \geq \lambda _n$ sont les valeurs propres du hessien $Hr$ d'une fonction de Busemann, alors $\lambda _{1}=0$, et $-b\leq \lambda _{i}\leq -a,$ pour $i\geq 2$. Ainsi, pour toute forme de degr\'e $k$, on a
\begin{eqnarray}
2HU(\alpha,\alpha )+\Delta U \vert \alpha \vert^{2} -|\nabla U |^2 |\alpha |^2 &\leq& C(2(1-k)a+(n-1)b-C)|\alpha |^2 \nonumber \\
&\leq& C(2a+(n-1)b-C) |\alpha |^2, \nonumber
\end{eqnarray}
puis on conclut comme avant.\cqfd

\Rq Comme on l'a d\'ej\`a sugg\'er\'e, on peut faire des raisonnements analogues si on consid\`ere le poids $e^{-U}$. Ainsi, si $M$ et $U$ v\'erifient les hypoth\`eses 1, 2 et 3 du th\'eor\`eme \ref{poids}, avec $1+2nC_1-nC_2>0$ (ou bien $C_1\leq 0$ et $1-nC_2>0$), alors la $L^2-$cohomologie \`a poids de $M$, associ\'ee au poids $e^{-U}$, est isomorphe \`a la cohomologie r\'eelle de $M$. De fa\c con similaire, si $M$ et $U$ sont comme dans le corollaire pr\'ec\'edent, avec $C>2nb-(n-1)a$, alors la $L^2-$cohomologie associ\'ee au poids $e^{-U}$ est aussi isomorphe \`a la cohomologie r\'eelle.


\section{Un d\'etour par les travaux d'Ohsawa}
Dans \cite{O} (voir aussi \cite{O-T}), T. Ohsawa d\'eveloppe une th\'eorie afin d'obtenir des th\'eor\`emes d'isomorphisme pour les groupes de cohomologie d'une large classe de vari\'et\'es complexes. Sa th\'eorie des "pseudo-paires de Runge" est en effet assez g\'en\'erale, et permet de comparer les groupes de $L^2$-cohomologie et de cohomologie ordinaire. Ici, nous nous contentons de voir que cette m\'ethode marche dans un cadre plus simple et l\'eg\`erement diff\'erent. Tout d'abord, nous avons besoin d'introduire la d\'efinition suivante, qui est l'\'equivalent de la notion de pseudo-paire de Runge chez Ohsawa :
\begin{definition}\label{PPR}
Soient $(M,g)$ une vari\'et\'e riemannienne compl\`ete, $k$ un entier, et $U : M\ra \R$ une fonction de classe $C^2$ d\'efinie sur $M$. On consid\`ere une suite de fonctions $U_p :M\ra \R$ d\'efinies sur $M$, de classe $C^2$. On dit que la suite $U_p$ est une bonne suite de poids par rapport \`a $U$ en degr\'e $k$ si elle satisfait aux conditions suivantes :
\be
 \item[(1)] Pour tout compact $L$ de $M$, il existe un entier $p_0$ (d\'ependant de $L$) tel que $U_p|L=U|L$ pour $p\geq p_0$ .
 \item[(2)](In\'egalit\'es de Poincar\'e \`a l'infini uniformes) Il existe un compact $K$ de $M$, et une constante $C_1>0$ tels que pour tout $p\geq 1$, l'on ait
  $$\forall \alpha \in C^\infty _0(\Lambda ^kT^*(M)),\, ||\alpha || ^2_{L^2_{U_p}}\leq C_1\left[ ||d\alpha ||^2_{L^2_{U_p}} +||\delta _{U_p}\alpha ||^2_{L^2_{U_p}} +\int _K \vert\alpha \vert ^2 dv_{U_p}\right] .$$
 \item[(3)] $L^2_{U_p}(\Lambda ^j T^*M) \subset L^2_{U_{p+1}}(\Lambda ^j T^*M)$,
  et il existe une constante $C_2>0$ telle que $||.||_{L^2_U(\Lambda ^j T^*M)}
  \leq C_2 ||.||_{L^2_{U_p}(\Lambda ^j T^*M)}$, pour tout $p\geq 1$, et $j=k-1,
  k$.  
\ee
\end{definition}
\Rq Il n'est pas difficile de voir que la condition (2) est \'equivalente \`a : il existe un compact $K'$ de $M$, et une constante $C'>0$ tels que pour tout $p\geq 1$, l'on ait
  $$\forall \alpha \in C^\infty _0(\Lambda ^kT^*(M\sm K')),\, ||\alpha || ^2_{L^2_{U_p}}\leq C'\left[ ||d\alpha ||^2_{L^2_{U_p}} +||\delta _{U_p}\alpha ||^2_{L^2_{U_p}}\right] ,$$ et dans la pratique, c'est cette condition qu'on v\'erifiera.

\begin{theoreme}[Ohsawa]\label{Ohsawa}
Soient $(M,g)$ une vari\'et\'e riemannienne compl\`ete, $k$ un entier, et $U : M\ra \R$ une fonction de classe $C^2$ d\'efinie sur $M$. On suppose qu'il existe une bonne suite de poids $U_p$ par rapport \`a $U$ en degr\'es $k$ et $k+1$ (c'est la m\^eme suite pour les deux degr\'es). Alors il existe un entier $p_0$ tel que les applications naturelles $H^k_{2, U_p}(M)\ra
H^k_{2, U}(M)$ soient des isomorphismes pour $p\geq p_0$.
\end{theoreme}

Comme notre cadre est un peu diff\'erent de celui consid\'er\'e par Ohsawa, et pour \^etre plus complet, nous refaisons la d\'emonstration de ce th\'eor\`eme en reprenant les arguments de la preuve de \cite[th\'eor\`eme 3.2, chapitre 2]{O}. Pour cela, nous avons besoin tout d'abord de quelques r\'esultats pr\'eliminaires. Dans la suite, $M$ est une vari\'et\'e riemannienne compl\`ete, et $U$ une fonction de classe $C^2$ d\'efinie sur $M$. En premier lieu, l'analogue de \cite[lemme 2.2, chapitre 2]{O} est 
\begin{prop}\label{marre}
On suppose qu'il existe une bonne suite de poids $U_p$ par rapport \`a $U$ en degr\'e $k$. Alors il existe un entier $p_0$ et une constante $C_3>0$ tels que pour tout $p\geq p_0$ l'on ait
$$\forall \alpha \in L^2_{U_p}(\Lambda ^kT^*M),\, \alpha \perp _{U} \mathcal{H}^k_U,\, ||\alpha || ^2_{L^2_{U_p}}\leq C_3\left[ ||d\alpha ||^2_{L^2_{U_p}} 
 +||\delta _{U_p}\alpha ||^2_{L^2_{U_p}} \right] .$$
\end{prop}
\dem Nous reprenons le raisonnement de \cite[lemme 2.2, chapitre 2]{O}, en supposant que la conclusion de la proposition est fausse. Alors il existe une suite de $k-$formes $\alpha _p$ v\'erifiant
$$\alpha _p\in L^2_{U_p},\,\, ||\alpha _p||_{L^2_{U_p}}=1,\,\, \alpha _p\perp \mathcal{H}^k_U,\,\, ||d\alpha _p||_{L^2_{U_p}}\leq 1/p,\,\, ||\delta _{U_p}\alpha _p||_{L^2_{U_p}}\leq 1/p.$$
Alors les $\alpha _p$ forment une suite born\'ee de $L^2_U$ d'apr\`es l'hypoth\`ese (3) (voir la d\'efinition \ref{PPR}), donc on peut en extraire une sous-suite qui converge faiblement dans $L^2_U$ vers une \'el\'ement $\alpha$ ; de plus, l'hypoth\`ese (1) montre que les $(d+\delta _U)(\alpha _p)$ sont born\'es sur les compacts, et donc on peut en particulier supposer, d'apr\`es le th\'eor\`eme de compacit\'e de Rellich, que les $\alpha _p$ convergent fortement vers $\alpha$ sur le compact $K$ (\`a une sous-suite pr\`es). Mais $\alpha$ est \`a la fois $(d+\delta _U)-$harmonique, et orthogonale \`a $\mathcal{H}^k_U$, donc $\alpha$ est nulle. Or, l'hypoth\`ese (2) entra\^\i ne
$$\frac{1}{C_1}-\frac{2}{p^2}\leq\int _K|\alpha _p |^2 dv_{U_p} ,$$
et en passant \`a la limite quand $p$ tend vers l'infini, on voit que
$$\frac{1}{C_1}\leq \int _K|\alpha |^2 dv_{U} ,$$ ce qui contredit le fait que $\alpha$ est nulle. \cqfd 

Avant d'\'enoncer la prochaine proposition, nous rappelons un r\'esultat d'H\"ormander \cite[th\'eor\`eme 1.1.4]{Ho} qui nous sera utile. Soient $H_i$, $i=1,2,3$ trois espaces de Hilbert, de normes $||.||_i$ respectivement. On consid\`ere deux op\'erateurs non born\'es et ferm\'es $T : H_1\ra H_2$ et $S : H_2\ra H_3$ de domaines respectifs $\rm{Dom}(T)$ et $\rm{Dom}(S)$ denses, et tels que $S\circ T=0$. On note par $T^*$ l'adjoint de $T$.
\begin{theoreme}[H\"ormander]\label{Hormander}
Soit $F$ un sous-espace ferm\'e de $H_2$ qui contient $\rm{Im}(T)$. On suppose qu'il existe une constante $c>0$ telle que
$$\forall \alpha \in \rm{Dom}(T^*)\cap\rm{Dom}(S)\cap F,\,\, ||\alpha ||^2_2\leq c\left[ ||T^*\alpha ||^2_1+||S\alpha ||^2_3\right] .$$
Alors si $\beta \in \rm{Im}(T^*)$, il existe $\alpha \in\rm{Dom}(T^*)$, avec $T^*\alpha =\beta$ et $||\alpha ||_2\leq c^{-1/2}||\beta ||_1.$
\end{theoreme}
Nous pouvons maintenant prouver l'analogue du th\'eor\`eme 2.3, chapitre 2 de \cite{O}, avec la m\^eme m\'ethode.
\begin{prop}
On suppose qu'il existe une bonne suite de poids $U_p$ par rapport \`a $U$ en degr\'e $k$. Soit $\gamma \in \rm{Ker}_{L^2_U}(d)$, de degr\'e $k-1$. Alors pour tout $\epsilon >0$, il existe un entier $p$ et un \'el\'ement $\tilde{\gamma}\in \rm{Ker}_{L^2_{U_{p}}}(d)$, de degr\'e $k-1$, tels que $$||\tilde{\gamma}-\gamma ||_{L^2_U}<\epsilon .$$
\end{prop}
\dem La conclusion de la proposition signifie que l'espace $\left( \cup _{p\geq 1}\rm{Ker}_{L^2_{U_{p}}}(d) \right) \cap L^2_U(\Lambda ^{k-1}T^*M)$ est dense dans $\rm{Ker}_{L^2_U}(d)\cap L^2_U(\Lambda ^{k-1}T^*M)$. Il suffit donc de prouver que l'orthogonal du premier espace dans le deuxi\`eme est nul : soit $\alpha\in \rm{Ker}_{L^2_U}(d)$ une forme de degr\'e $k-1$ telle que
$$\forall p\geq 1,\,\,\forall \phi \in \rm{Ker}_{L^2_{U_{p}}}(d),\,\, \lan \alpha ,\phi\ran _{L^2_U}=0.$$ On va montrer que $\alpha =0$, et pour cela on va montrer qu'il existe une forme $\beta$ de degr\'e $k$ dans $L^2_U$ telle que $\alpha =\delta _U \beta$.

D'abord, l'hypoth\`ese (3) montre que pour tout $p$, la forme lin\'eaire $\lan \alpha ,.\ran _{L^2_U}$ est continue sur $L^2_{U_p}$, de norme inf\'erieure ou \'egale \`a $C_2 ||\alpha ||_{L^2_U}$, donc on en d\'eduit l'existence de $\alpha _p\in L^2_{U_p}$ telle que
$$\lan \alpha , .\ran _{L^2_U}=\lan \alpha _p , .\ran _{L^2_{U_p}},\,\, \mathrm{et} \,\, ||\alpha _p ||_{L^2_{U_p}}\leq C_2 ||\alpha ||_{L^2_U}.$$
Plus simplement, on a $\alpha _p=e^{U-U_p}\alpha$. D'apr\`es l'hypoth\`ese (1), on voit alors que
$$\forall \phi \in C^\infty _0(\Lambda ^kT^*M),\,\, \lan \alpha _p, \phi\ran _{L^2_U}\ra \lan \alpha, \phi\ran _{L^2_U} ,$$ puisque pour tout compact fix\'e, il y a en fait \'egalit\'e entre $\alpha$ et $\alpha _p$  pour tous les $p$ assez grands.

D'autre part, l'hypoth\`ese d'orthogonalit\'e faite sur $\alpha$, et la d\'efinition des $\alpha _p$ montrent que $\alpha _p$ est orthogonale \`a $\rm{Ker}_{L^2_{U_{p}}}(d)$ dans $L^2_{U_p}$, donc est dans l'adh\'erence de l'image de $\delta _{U_p}$. Or, cette image est ferm\'ee gr\^ace aux in\'egalit\'es de Poincar\'e \`a l'infini (2), par un raisonnement proche de celui de la preuve de la proposition \ref{exact} : supposons en effet qu'il existe une suite $\psi _l$ de $k-$formes de $L^2_{U_p}$ telle que la suite $\delta _{U_p}\psi _l$ soit dans $L^2_{U_p}$, et converge dans cet espace vers une $(k-1)-$forme $\omega$. L'in\'egalit\'e de Poincar\'e \`a l'infini montre que $0$ n'est pas dans le spectre essentiel du laplacien $\Delta _{U_p,k}$ agissant sur les $k-$formes (o\`u $\Delta _{U_p}=(d+\delta _{U_p})^2$), donc on a un op\'erateur de Green $G$ comme dans la proposition \ref{Green}. Quitte \`a remplacer $\psi _l$ par $dG_{-}\psi _l$, on peut supposer que les $\psi _l$ sont exactes. On utilise alors encore une fois l'hypoth\`ese sur le spectre essentiel pour en d\'eduire l'existence d'une constante $C>0$ telle que 
$$\forall \psi \in \rm{Dom}((d+\delta _{U_p})_k), \, \psi \bot \rm{Ker}(\Delta _{U_p,k}),\,  C\Vert \psi \Vert ^2 _{L^2_{U_p}}\leq  \Vert d\psi \Vert ^2 _{L^2_{U_p}} +\Vert \delta _{U_p} \psi \Vert ^2 _{L^2_{U_p}}.$$
On en d\'eduit facilement que $\psi _l$ est une suite de Cauchy dans $L^2_{U_p}$, et qu'elle converge vers un \'el\'ement $\psi$, avec $\delta _{U_p}\psi =\omega$. Ceci finit de montrer que l'image de $\delta _{U_p}$ est ferm\'ee dans $L^2_{U_p}$. On a donc $$\alpha _p\in \rm{Im}(\delta _{U_p}).$$
On va appliquer le th\'eor\`eme d'H\"ormander \ref{Hormander} avec $H_i =L^2_{U_p}(\Lambda ^{k+i-2}T^*M)$, $i=1,2,3$, $T=S=d$, et $F$ l'ensemble des $k$-formes de $L^2_{U_p}$ qui sont orthogonales \`a $\mathcal{H }^k_U$ dans $L^2_U(\Lambda ^kT^*M)$. $F$ est clairement ferm\'e, et de plus il contient l'image de $T=d$ : en effet, soient $\xi$ un \'el\'ement de $\mathcal{H }^k_U$, et $\gamma =T\eta =d\eta \in L^2_{U_p}(\Lambda ^kT^*M)\subset L^2_U(\Lambda ^kT^*M) $ un \'el\'ement de $\rm{Im}(T)$, avec $\eta \in L^2_{U_p}(\Lambda ^{k-1}T^*M)\subset L^2_U(\Lambda ^{k-1}T^*M)$. Alors on a :
\begin{eqnarray*}
\lan \xi ,\gamma \ran _{L^2_U} &=&\lan \xi ,d\eta \ran _{L^2_U}\\
                               &=&\lan \delta _U\xi ,\eta \ran _{L^2_U}=0,
\end{eqnarray*}
l'int\'egration par parties qui fait passer de la premi\`ere ligne \`a la deuxi\`eme \'etant justifi\'ee car $(M,g)$ est compl\`ete. On a donc bien $\rm{Im}(T)\subset F$. Enfin, l'estimation dont on a besoin pour appliquer le th\'eor\`eme d'H\" ormander est obtenue dans la proposition \ref{marre}, avec la m\^eme constante $c=C_3$ pour tous les $p$ assez grands. Ceci montre l'existence d'une constante $C_4$ et, pour tout $p$, d'une forme $\beta _p$ de $L^2_{U_p}$ telles que
$$\delta _{U_p}\beta _p=\alpha _p,\,\, \mathrm{et} \,\, ||\beta _p||_{L^2 _{U_p}}\leq C_4 ||\alpha _p||_{L^2 _{U_p}}.$$
Il est alors facile de voir que les $\beta _p$ sont born\'ees dans $L^2_U$ (par $C_2^2C_4||\alpha ||_{L^2 _U}$), donc que les $\beta _p$ convergent (\`a une sous-suite pr\`es, not\'ee encore $\beta _p$) faiblement vers  une forme $\beta$ dans $L^2_U$. De plus, on a $\delta _U\beta =\alpha$ : en effet, pour toute forme $\phi$ lisse \`a support compact, et pour tout $p$ assez grand (d\'ependant du support de $\phi$), on a d'apr\`es (1)
$$\lan \beta _p,d\phi\ran _{L^2_{U_p}}=\lan \beta _p,d\phi\ran _{L^2_U}.$$
Le membre de droite de cette \'egalit\'e tend vers $\lan \beta ,d\phi\ran _{L^2_U}$ par convergence faible. Quant au membre de gauche, on a 
$$\lan \beta _p,d\phi\ran _{L^2_{U_p}}=\lan \delta _{U_p}\beta,\phi\ran _{L^2_{U_p}}=\lan \alpha _p,\phi \ran _{L^2_{U_p}}=\lan\alpha ,\phi \ran _{L^2_U}.$$
Ceci ach\`eve la preuve de la proposition. \cqfd 

\emph{D\'emonstration du th\'eor\`eme \ref{Ohsawa}}. Gr\^ace \`a (1) et (2), l'op\'erateur $(d+\delta _U)$ v\'erifie aussi une in\'egalit\'e de Poincar\'e \`a l'infini, donc $H^k_{2,U}(M)$ est de dimension finie. Il suffit donc de montrer que pour $p$ assez grand, l'application $H^k_{2,U_p}(M)\ra H^k_{2,U}(M)$ est injective et d'image dense, et ces deux faits d\'ecoulent des deux propositions pr\'ec\'edentes.\cqfd 

\section{$L^2$-cohomologie des vari\'et\'es de volume fini, \`a courbure n\'egative pinc\'ee}

 \subsection{Th\'eor\`eme principal}
 Dans cette partie, $(M^{n},g)$ est une vari\'et\'e riemannienne compl\`ete de dimension $n$, de volume fini, \`a courbure sectionnelle $K$ n\'egative pinc\'ee : il existe des constantes $a$ et $b$, telles que l'on ait $-b^{2}\leq K\leq -a^{2} <0$. Dans la troisi\`eme partie (corollaire \ref{poidsfini}), on a donn\'e une interpr\'etation topologique des espaces de $L^2$-cohomologie \`a poids de ces vari\'et\'es, si le poids est bien choisi. Ici, nous allons comparer les espaces de $L^2$-cohomologie aux espaces de $L^2$-cohomologie \`a poids, gr\^ace aux r\'esultats de la partie pr\'ec\'edente. Des travaux d'Ohsawa, on d\'eduit d'abord la proposition suivante :
\begin{prop}\label{amelioration}
Soit $(M^{n},g)$, une vari\'et\'e riemannienne compl\`ete de dimension $n$, de volume fini, \`a courbure sectionnelle $K$ n\'egative pinc\'ee : il existe des constantes $a$ et $b$, telles que l'on ait $-b^{2}\leq K\leq -a^{2} <0$. On suppose que pour un entier $k<(n-1)/2$, on a  $(n-1-k)a-kb=:\varepsilon>0.$ Alors on a l'isomorphisme:
$$H^{n-k}_2(M)\simeq H^{n-k}_c(M).$$
\end{prop}
\dem On va appliquer le th\'eor\`eme \ref{Ohsawa} avec une suite de poids bien choisis. Fixons d'abord une constante $C>0$ qu'on choisira ensuite assez grande, et consid\'erons une fonction $F :\R\to \R$ de classe $C^2$ qui v\'erifie
$$F(t)=0 \quad\mathrm{si\,\, t\leq 0},$$
$$F(t)=Ct\quad\mathrm{si\,\, t\geq 1},$$
$$F'\geq 0.$$
 On peut, sans perdre de g\'en\'eralit\'es, supposer que $M$ a un seul bout, et on note par $r$ une fonction de Busemann associ\'ee. On consid\`ere alors une suite $U_p : M\ra \R$ de fonctions de classe $C^2$ d\'efinies en $x\in M$ par $$U_p(x)=F(r(x)-p).$$
Si $x$ est fix\'e dans $M$, pour $p$ assez grand, on a $U_p(x)=0$, c'est-\`a-dire que l'hypoth\`ese (1) de la d\'efinition \ref{PPR} est v\'erifi\'ee par $U\equiv 0$ et la suite $U_p$. L'hypoth\`ese (3) est aussi facile \`a v\'erifier avec $C_2=1$. Nous allons voir que la condition (2) est aussi satisfaite en degr\'es $n-k$ et $n-k+1$. En effet, d'apr\`es la formule d'int\'egration par parties du corollaire \ref{IPPpoids}, en prenant $f=r$, si $\alpha$ est une forme \`a support compact contenu dans le bout, on a 
\begin{multline} \label{IPP4}
\int _{M} \left[ 2Hr(\alpha,\alpha )+\Delta r \vert \alpha \vert^{2} -F'(r-p) |\alpha |^2 \right] \rm{d}v_{U_p} = \\ 2\int_{M} \left[\langle i_{\nabla r}d\alpha ,\alpha \rangle +\langle i_{\nabla r}\alpha ,\delta _{U_p} \alpha \rangle \right ] \rm{d}v_{U_p}.  
\end{multline}
Notons $\lambda _1\geq \ldots\geq\lambda _n$ les valeurs propres du hessien $Hr$. D'apr\`es l'hypoth\`ese sur les courbures sectionnelles, le th\'eor\`eme de comparaison des hessiens de Greene et Wu \cite{G-W} nous donne $\lambda _1=0$, et $-b\leq\lambda _i\leq -a$, $i\geq 2$ sur le bout. Ainsi, en raisonnant comme dans \cite[corollaire 5.4]{B-B}, si $\alpha$ est de degr\'e $n-k$, on a
$$2Hr(\alpha,\alpha )+\Delta r \vert \alpha \vert^{2}\leq -\left[ (n-k-1)a-kb\right] |\alpha |^2=-\epsilon |\alpha |^2.$$
Cette in\'egalit\'e reste \textit{a fortiori} vraie si $\alpha$ est de degr\'e $n-k+1$. Comme on a $F'\geq 0$, on d\'eduit de la formule d'int\'egration par parties \ref{IPP4} que pour toute $(n-k)$ ou $(n-k+1)$-forme $\alpha$ \`a support dans le bout :
$$\sqrt \varepsilon  \Vert \alpha \Vert _{L^2_{U_p}} \leq 2 \left[ \Vert d\alpha \Vert _{L^2_{U_p}}+\Vert \delta _{U_p} \alpha \Vert _{L^2_{U_p}}\right].$$ 
Ceci montre que l'hypoth\`ese (2) de la d\'efinition \ref{PPR} est aussi v\'erifi\'ee pour les degr\'es $n-k$ et $n-k+1$. Par cons\'equent, la suite $U_p$ est une bonne suite de poids pour $U\equiv 0$ en degr\'es $n-k$ et $n-k+1$, donc d'apr\`es le th\'eor\`eme d'Ohsawa \ref{Ohsawa},  pour $p$ assez grand, on a $$H^{n-k}_{2, U_p}(M)\simeq H^{n-k}_2(M).$$
Or, si on est assez loin sur le bout, la fonction $U_p$ est \'egale \`a $Cr$, et par le corollaire \ref{poidsfini}, on sait dans ce cas que si la constante $C$ est assez grande, alors $$H^*_{2, U_p}(M)\simeq H^*_c(M).$$ Ceci ach\`eve la preuve de la proposition.\cqfd

\Rq L'hypoth\`ese $(n-1-k)a-kb>0$ est exactement celle qui permet de savoir que $0$ n'est pas dans le spectre essentiel du laplacien $\Delta _k$ (voir \cite[th\'eor\`eme 4.2]{Do-X} et \cite[corollaire 5.4]{B-B}). On verra que le r\'esultat de la proposition \ref{amelioration} est optimal par rapport au pincement.

Nous pouvons maintenant \'enoncer notre r\'esultat principal.
\begin{theoreme}\label{cohomologie}
Soit $(M^{n},g)$, une vari\'et\'e riemannienne compl\`ete de dimension $n$, de volume fini, \`a courbure sectionnelle $K$ n\'egative pinc\'ee : il existe des constantes $a$ et $b$, telles que l'on ait $-b^{2}\leq K\leq -a^{2} <0$. On suppose que $na-(n-2)b>0,$ alors on a les isomorphismes entre espaces vectoriels de dimension finie :
$$H^k_2(M)\simeq \left\{ \begin{array}{lll}
     H^k(M), & \textrm{si $k<(n-1)/2$,}\\
{\rm Im} (H^{n/2}_c(M)\ra H^{n/2}(M)), & \textrm{si $k=n/2$,}\\
     H^k_c(M), & \textrm{si $k> (n+1)/2.$}
     \end{array} \right. $$
Si de plus la courbure est constante ($a=b$), et la dimension est impaire, alors on a $H^{(n\pm 1)/2}_2(M)\simeq {\rm Im} (H^{(n\pm 1)/2}_c(M)\ra H^{(n\pm 1)/2}(M))$ 
\end{theoreme}
\dem
On peut tout d'abord supposer que $M$ est orientable. Sinon, il suffit de consid\'erer le rev\^etement \`a deux feuillets $\tilde M$ qui, lui, est orientable, en remarquant que si $\tau$ est l'automorphisme du rev\^etement, alors la $L^2-$cohomologie de $M$ est isomorphe \`a la $L^2-$cohomologie $\tau -$invariante de $\tilde M$. L'avantage qu'on en tire est que sur une vari\'et\'e orientable, l'op\'erateur $*$ de Hodge r\'ealise une isom\'etrie entre les $k$ et les $(n-k)-$formes harmoniques. Ainsi, on peut se limiter aux degr\'es $k$ sup\'erieurs ou \'egaux \`a $n/2$. 
\be
 \item Pour les degr\'es $k>(n+1)/2$, la conclusion du th\'eor\`eme d\'ecoule de la proposition \ref{amelioration}
 \item Supposons que $n=2m$, et traitons le cas $k=n/2=m$. On sait que, sous nos hypoth\`eses de courbure, $0$ n'est pas dans le spectre essentiel du laplacien agissant sur les formes \cite{B-B}. Les deux suites exactes de la premi\`ere partie (corollaire \ref{se}) sont donc vraies pour tous les degr\'es, et en particulier, pour tout ouvert born\'e $D$, on a :
$$H^m_c(D) \ra H^m_2(M) \ra H^m_2(M\sm D),$$
et
$$H^m_2(M\sm D,\partial D) \ra H^m_2(M)\ra H^m(D).$$ 
Sans perdre de g\'en\'eralit\'es, on peut supposer que $M$ a un seul bout, et on note par $r$ une fonction de Busemann associ\'ee. On prend alors $D=\{ r<c\}$ pour un r\'eel $c$ fix\'e, et on va montrer que  $H^m_2(M\sm D)=H^m_2(M\sm D,\partial D)=0,$ si on choisit $c$ assez grand. Alors les deux suites exactes ci-dessus donneront : $H^{m}_{c}(M)$ se surjecte dans $H^m_2(M)$, et $H^m_2(M)$ s'injecte dans $H^{m}(M)$. Mais d'apr\`es Anderson \cite{An} (voir aussi \cite{L1} ou \cite{C4}), ${\rm Im} (H^{m}_{c}(M)\ra H^{m}(M))$ s'injecte toujours dans $H^m_2(M)$, donc on en d\'eduira que $H^m_2(M)\simeq {\rm Im} (H^{m}_{c}(M)\ra H^{m}(M)).$

Il nous reste donc \`a montrer que $H^m_2(M\sm D)=H^m_2(M\sm D,\partial D)=0.$ Pour cela, on consid\`ere la vari\'et\'e double $X_c=(M\sm D)\sharp (-M\sm D)$ obtenue en recollant  deux copies de $M\sm D$ le long de $\partial D$. $X_c$ porte une m\'etrique naturelle qui co\" incide avec celle de $M\sm D$ sur chacune des copies. Cette m\'etrique n'est pas lisse, mais seulement lipschitzienne ; cependant, ceci est suffisant pour d\'efinir l'op\'erateur $\delta$ et toute l'analyse pr\'ec\'edente reste valable. La $L^2-$cohomologie de $X_c$ v\'erifie alors (voir par exemple \cite{C1})
$$H^{*}_2(X_c)=H^{*}_2(M\sm D)\oplus H^{*}_2(M\sm D,\partial D).$$
Il nous suffit donc de montrer que $H^m_2(X_c)=0$. D'apr\`es la premi\`ere partie de la d\'emonstration, on sait d\'ej\`a que 
  $$H^k_2(X_c)\simeq \left\{ \begin{array}{ll}
     H^k(X_c), & \textrm{si $k<m$,}\\
     H^k_c(X_c), & \textrm{si $k>m.$}
     \end{array} \right. $$
$X_c$ est diff\'eomorphe \`a un produit $\R \times \partial D$, donc ses nombres de Betti $b_k(X_c)$ sont les m\^emes que ceux de $\partial D$. La caract\'eristique d'Euler $L^2$ de $X_c$ est par cons\'equent donn\'ee par 
$$\chi _{2}(X_c)=2\sum _{k<m}(-1)^k b_k(\partial D )+(-1)^m\rm{dim}(H^m_2(X_c)),$$ o\`u $b_k(\partial D )$ est le $k-$\`eme nombre de Betti de $\partial D$. Or, d'apr\`es \cite[th\'eor\`eme C]{B-B}, on a :
$$\chi _2(X_c)=\int _{X_c} \omega +2\int _{\partial D}P(II)+2\sum _{k<m}(-1)^k b_k(\partial D ),$$ o\`u $\omega$ est la forme d'Euler, et $P(II)$ est un polyn\^ome en la courbure et en la seconde forme fondamentale de $\partial D$. On en d\'eduit donc que :
$$(-1)^m\rm{dim}(H^m_2(X_c))=\int _{X_c} \omega +2\int _{\partial D}P(II) .$$
Mais le terme de droite de cette \'egalit\'e tend vers z\'ero lorsque $c$ tend vers l'infini, donc le terme de gauche, qui est un entier, est nul pour $c$ assez grand ; on obtient bien l'annulation souhait\'ee.
 \item Il ne nous reste plus qu'\`a consid\'erer le cas o\`u la courbure est constante (par exemple \'egale \`a $-1$), $n=2m+1$ est impair, $k=(n-1)/2=m$. On peut toujours supposer qu'on a un seul bout, et donc qu'en dehors d'un ouvert $D$, $M\sm D$ est de la forme $]0,\infty [ \times \partial D$, portant la m\'etrique $$dr^2+e^{-2r}d\theta ^2,$$
o\`u $r$ est une fonction de Busemann, et $d\theta ^2$ est une m\'etrique sur la vari\'et\'e $\partial D$. D'apr\`es \cite{Do-X}, $0$ n'est pas dans le spectre essentiel du laplacien en degr\'e $m-1$. Donc d'apr\`es les r\'esultats de la premi\`ere partie, on a les deux suites exactes
 $$H^{m}_{c}(D) \ra H^m_2(M) \ra H^m_2(M\sm D),$$
 $$H^m_2(M\sm D,\partial D) \ra H^m_2(M)\ra H^{m}(D).$$ 
De nouveau, il nous suffit de montrer l'annulation de la $L^2-$cohomologie \`a l'infini. On va proc\'eder comme dans la d\'emonstration de \cite[lemme 5.3]{C4}. D'abord, $]0,\infty [ \times \partial D $ est finiment recouvert par $ ]0,\infty [ \times \T ^{n-1}$, avec $\T ^{n-1}$ un tore plat. Il suffit donc de montrer le r\'esultat pour $]0,\infty [ \times \T ^{n-1}$. Comme le tore $\T ^{n-1}$ agit par isom\'etries sur $]0,\infty [ \times \T ^{n-1}$, et que cette action engendre des vecteurs de Killing de longueur born\'ee, \cite[th\'eor\`eme 3]{Hi} montre que toute forme harmonique $L^2$ v\'erifiant la condition absolue ou relative est invariante sous cette action. Ainsi, si $\alpha$ est une telle $m-$forme, on peut \'ecrire $\alpha = \beta (r)+dr \wedge \gamma (r)$, o\`u $\beta (r)$ et $\gamma (r)$ sont vues comme des formes d\'efinies sur $\T ^{n-1}$, d\'ependant du param\`etre $r$, et on a :
$$d\alpha = dr\wedge \frac{\partial \beta }{\partial r}=0,$$
$$\delta \alpha = -\frac{\partial \gamma}{\partial r}+2\gamma =0,$$
donc si $\alpha$ v\'erifie la condition absolue $\gamma (0)=0$, on obtient finalement $\alpha = \beta (0).$ Mais cette forme n'est de carr\'e int\'egrable que si elle est nulle. De m\^eme, si $\alpha$ v\'erifie la condition relative, $\alpha$ est forc\'ement nulle.    
\ee\cqfd

{\bf Remarques.}
\be
\item[i)] Pour les vari\'et\'es hyperboliques de volume fini, on retrouve les r\'esultats de Zucker \cite{Z}, et de Mazzeo et Phillips \cite{M-P} concernant la $L^2-$cohomologie (ces deux derniers auteurs traitent le cas plus g\'en\'eral des vari\'et\'es hyperboliques g\'eom\'etriquement finies).
\item[ii)] Les r\'esultats pr\'ec\'edents, et tous ceux qui suivent, restent vrais si on perturbe la m\'etrique sur un compact.
\item[iii)] Pour avoir la conclusion seulement pour les degr\'es inf\'erieurs ou \'egaux \`a un $k<(n-1)/2$, il suffit d'avoir $(n-k-1)a-kb>0$ (voir la proposition \ref{amelioration}. Mieux encore, il suffit d'imposer ces conditions de pincement uniquement sur les courbures radiales, car ce sont celles qui interviennent dans le th\'eor\`eme de comparaison des hessiens de Greene et Wu.
\item[iv)] Les degr\'es $k=(n\pm 1)/2$ quand $n=2m+1$ est impair et la courbure variable ne sont pas pris en compte : on verra en effet \`a la fin de cette partie un exemple montrant qu'on ne peut pas esp\'erer de r\'esultats analogues pour ces degr\'es d\`es que la courbure varie. 
\item[v)] Si la courbure est constante, et si la dimension est paire, ou bien si on est en dimension deux, l'annulation de la $L^2-$cohomologie \`a l'infini en degr\'e moiti\'e s'obtient de mani\`ere plus simple, par invariance conforme de la $L^2-$cohomologie en ce degr\'e : dans ce cas, chaque bout est conform\'ement \'equivalent \`a un cylindre, et il est bien connu que la $L^2-$cohomologie d'un cylindre est nulle.   
\ee

  \subsection{Etude d'un exemple : le produit doublement tordu}
On va maintenant \'etudier une famille d'exemples qui montre que la conclusion du th\'eor\`eme est fausse, d'une part pour les degr\'es proches du degr\'e moiti\'e si la dimension est impaire, d\`es que la courbure varie, et d'autre part pour tous les degr\'es (autres que $0$ et $n$) si la courbure n'est pas assez pinc\'ee. Consid\'erons le tore $\T ^i$ de dimension $i$, notons $h_i$ la m\'etrique induite par celle de $\R ^i$. Soient $n\geq 3$ un entier, $1\leq k<n-1$ un autre entier, et $0<a< b$ deux r\'eels. On munit la vari\'et\'e $\Omega =]0,\infty [\times \T ^{n-k-1}\times\T ^k$ de la m\'etrique $$g=dr^2+e^{-2ar}h_{n-k-1}+e^{-2br}h_{k}.$$
Ballmann et Br\"uning ont consid\'er\'e cet exemple \cite[exemples 3.19, 4.2 et 5.5]{B-B}. Ils montrent notamment que la courbure sectionnelle $K$ v\'erifie $-b^2\leq K\leq -a^2,$ et que $0$ est dans le spectre essentiel du laplacien $\Delta _k$ si et seulement si $(n-k-1)a-kb=0$.

Nous allons d\'eterminer l'espace des $p-$formes harmoniques sur $(\Omega ,g)$, en raisonnant comme dans \cite[lemme 5.3]{C4} et comme \`a la fin de la d\'emonstration du th\'eor\`eme pr\'ec\'edent. Tout d'abord, si $\alpha$ est une $p-$forme harmonique v\'erifiant la condition absolue ou relative sur le bord $\{ 0\}\times \T ^{n-k-1}\times\T ^k$, alors d'apr\`es le r\'esultat de Hitchin \cite[th\'eor\`eme 3]{Hi}, $\alpha$ est invariante sous l'action de $\T ^{n-k-1}\times\T ^k$. Si $x_1,\ldots ,x_{n-k-1},y_1,\ldots ,y_k$ est un syst\`eme de coordonn\'ees sur $\T ^{n-1}=\T ^{n-k-1}\times\T ^k$, on peut \'ecrire  
$$\alpha =\sum _{I,J}a_{I,J}(r)dx^I\wedge dy^J+dr\wedge \sum _{I,J}b_{I,J}(r)dx^I\wedge dy^J,$$
o\`u $I=(i_1,\ldots ,i_l)$ et $J=(j_1,\ldots ,j_m)$ sont des multi-indices, avec $l+m=p$ pour la premi\`ere somme et $l+m=p-1$ pour la deuxi\`eme, $l\leq n-k-1$, $m\leq k$, $dx^I\wedge dy^J =dx^{i_1}\wedge \ldots dx^{i_l}\wedge dy^{j_1}\wedge \ldots \wedge dy^{j_m}$, et enfin o\`u $a_{I,J}(r)$ , $b_{I,J}(r)$ sont des fonctions qui ne d\'ependent que de $r$. On peut r\'e\'ecrire cette formule sous la forme :
$$\alpha =\sum _I \beta _I(r)+\sum _I dr\wedge \gamma _I(r),$$ o\`u chaque multi-indice $I=(i_1,\ldots ,i_{l_I},j_1,\ldots ,j_{m_I})$ est de longueur $l_I+m_I=p$ (respectivement $l_I+m_I=p-1$) dans la premi\`ere somme (respectivement dans la deuxi\`eme somme). On a alors :
$$0=d\alpha = dr \wedge \sum _I\frac{\partial \beta _I}{\partial r},$$
$$0=\delta \alpha =-\sum _I [\frac{\partial \gamma _I}{\partial r}+(2(l_Ia+m_Ib)-(n-k-1)a-kb)\gamma _I].$$
Si $\alpha$ v\'erifie la condition absolue au bord, on a donc forc\'ement $$\alpha =\sum _I\beta _I,$$ avec les $\beta _I$ ind\'ependantes de $r$. Chaque $\beta _I$ \'etant de carr\'e int\'egrable si et seulement si $2(l_Ia+m_Ib)<(n-k-1)a+kb,$ on en d\'eduit finalement :
\begin{eqnarray}\label{dimension}
\rm{dim}(H^p_2(\Omega )) &=& \sharp \{I=(i_1,\ldots ,i_l,j_1,\ldots ,j_m)/\, l\leq n-k-1,\, m\leq k,\nonumber \\ 
 & &\quad l+m=p \quad \mathrm{et}\quad 2(la+mb)<(n-k-1)a+kb\},\nonumber\\
\end{eqnarray}
et comme  par dualit\'e $H^p_2(\Omega ,\partial \Omega )\simeq H^{n-p}_2(\Omega )$, 
\begin{eqnarray} 
\rm{dim}(H^p_2(\Omega ,\partial \Omega )) &=& \sharp \{I=(i_1,\ldots ,i_l,j_1,\ldots ,j_m)/\, l\leq n-k-1,\, m\leq k,\nonumber \\ 
 & &\quad l+m=n-p \quad \mathrm{et}\quad 2(la+mb)<(n-k-1)a+kb\}.\nonumber
\end{eqnarray}
On peut alors voir qu'en dimension impaire $n=2s+1$, la conclusion du th\'eor\`eme \ref{cohomologie} concernant les vari\'et\'es hyperboliques de volume fini n'est plus v\'erifi\'ee en courbure variable pour les degr\'es $s$ et $s+1$. Consid\'erons en effet le double produit tordu $\Omega =]0,\infty [\times \T ^s\times\T ^s$ muni de la m\'etrique $g=dr^2+e^{-2ar}h_{s}+e^{-2br}h_{s}.$ L'\'equation \ref{dimension} qui donne la dimension de $H^s_2(\Omega )$ montre que cette dimension n'est jamais nulle, d\`es que $b$ est strictement plus grand que $a$: le multi-indice qui a des $1$ pour $s$ premi\`eres composantes et des $0$ ailleurs est par exemple toujours solution dans ce cas. Ainsi, sur la vari\'et\'e double $\Omega \sharp (-\Omega)$, de dimension $2s+1$, il y a au moins une $s-$forme harmonique en degr\'e $s$, car $H^*_2(\Omega \sharp (-\Omega))=H^*_2(\Omega)\oplus H^*_2(\Omega ,\partial \Omega).$ Mais l'image de la cohomologie relative dans la cohomologie est nulle.

Ensuite, nous nous contentons, pour simplifier, de consid\'erer le cas $n=4$. Le cas g\'en\'eral est similaire, et n'apporte rien de plus. Avec les notations pr\'ec\'edentes, nous choisissons $k=1$, et nous pouvons supposer que $b=1$. L'\'egalit\'e \ref{dimension} montre que la $L^2-$cohomologie absolue de degr\'e un $H^1_2(\Omega )$ est de dimension 2 ou 3 selon les valeurs de $a$ : les formes $dx_1$ et $dx_2$ sont toujours dans $H^1_2(\Omega )$ et non nulles, alors que $dy_1$ n'est solution que si $a>1/2$. Quant \`a la $L^2-$cohomologie relative, elle est toujours nulle. Comme la cohomologie usuelle de degr\'e $1$ de $\Omega \sharp (-\Omega)=\R \times \T ^3$ est de dimension 3, l'isomorphisme $H^1(\Omega \sharp (-\Omega))\simeq H^1_2(\Omega \sharp (-\Omega))$ n'est vrai que si $a>1/2$. Remarquons que $a=1/2$ est exactement la valeur pour laquelle $0$ est dans le spectre essentiel de $\Delta _1$. De m\^eme, on voit qu'en degr\'e $2$, l'isomorphisme entre la $L^2-$cohomologie et l'image de la cohomologie \`a support compact dans la cohomologie n'est vrai que si $a>1/2$.

 \section{Quelques autres g\'eom\'etries}

  \subsection{Certains produits tordus}

On consid\`ere ici des vari\'et\'es dont chacun des bouts est isom\'etrique \`a une vari\'et\'e de la forme $B=]0,\infty[\times K$ (avec $K$ une vari\'et\'e compacte sans bord) munie d'une m\'etrique qui s'\'ecrit 
$$dr^{2}+F(r)^{2}d\theta ^{2},$$ o\`u $F$ est une fonction strictement positive, et $d\theta ^{2}$ est une m\'etrique sur $K$. Les formes harmoniques de telles vari\'et\'es ont par exemple \'et\'e \'etudi\'ees dans \cite{D}, \cite{Br} et \cite{C2}, o\`u les auteurs obtiennent des r\'esultats d'annulation et de finitude.

Nous allons donner une interpr\'etation topologique de ces espaces dans certains cas particuliers. Dans un rep\`ere orthonormal local $e_{1}=\nabla r,e_{2},\ldots ,e_{n}$ de $B$ (avec $e_{2},\ldots ,e_{n}$ un rep\`ere orthogonal de $K$ pour la m\'etrique $d\theta ^{2}$), le hessien de la fonction $r$ se diagonalise : il vaut z\'ero sur $e_{1}=\nabla r$, et $F'/F$ sur les autres vecteurs. Ainsi, si $\alpha$ est une $k-$forme, on a :
$$2Hr(\alpha ,\alpha )+\Delta r\vert \alpha \vert ^{2}=\frac{F'}{F}[(2k-(n-1))\vert \alpha \vert ^{2}-2\vert i_{\nabla r}\alpha \vert ^{2}].$$
Si on suppose que $F'/F\geq C>0$ pour une constante $C$, alors en utilisant la double in\'egalit\'e $0\leq \vert i_{\nabla r}\alpha \vert ^{2}\leq \vert \alpha \vert ^{2}$, on obtient 
\begin{eqnarray*}
2Hr(\alpha ,\alpha )+\Delta r\vert \alpha \vert ^{2}&\geq &C[2k-2-(n-1)]\vert \alpha \vert ^{2},\\
2Hr(\alpha ,\alpha )+\Delta r \vert \alpha \vert ^{2}&\leq &C[2k-(n-1)]\vert \alpha \vert ^{2}.
\end{eqnarray*}
On peut aussi faire un raisonnement analogue si on suppose que $F'/F\leq -C<0$. Ainsi, les m\'ethodes des paragraphes pr\'ec\'edents permettent d'affirmer :
\begin{theoreme}\label{produit tordu}
Soit $(M^{n},g)$ une vari\'et\'e riemannienne. On suppose que $M$ se d\'ecompose en $D\cup _{i} B_{i}$ o\`u $D$ est une vari\'et\'e compacte \`a bord, et chaque bout $B_{i}$ est isom\'etrique \`a un produit tordu de la forme $]0,\infty [\times K _{i}$ portant la m\'etrique $dr^{2}+F_{i}^{2}(r)d\theta _{i}^{2}$.
\bd
 \item[a)] S'il existe une constante $C>0$ telle que pour tout $i$ on ait $F'_{i}/F_{i}>C$, alors
$$H^k_2(M)\simeq \left\{ \begin{array}{ll}
     H^k_c(M), & \textrm{si $k\leq (n-1)/2$,}\\
     H^k(M), & \textrm{si $k\geq (n+1)/2.$}
     \end{array} \right. $$
 \item[b)] S'il existe une constante $C>0$ telle que pour tout $i$ on ait $F'_{i}/F_{i}<-C$, alors
$$H^k_2(M)\simeq \left\{ \begin{array}{lll}
     H^k(M), & \textrm{si $k< (n-1)/2$,}\\
{\rm Im} (H^k_c(M)\ra H^k(M)), & \textrm{si $k=n/2$}\\
     H^k_c(M), & \textrm{si $k> (n+1)/2.$}
     \end{array} \right. $$
\ed
\end{theoreme}
\dem On fait la d\'emonstration par exemple uniquement pour le cas b). Par dualit\'e, il suffit de consid\'erer les degr\'es $k\leq n/2$. Pour $k<(n-3)/2$, on raisonne comme dans la preuve du th\'eor\`eme \ref{poids} : les estim\'ees sur le hessien des fonctions $F_i$ montrent qu'on a une in\'egalit\'e de Poincar\'e \`a l'infini pour $k<(n-1)/2$, donc on a une suite exacte comme dans le corollaire \ref{se}. Les m\^emes estim\'ees montrent que pour $k<(n-1)/2$, la $L^2$-cohomologie relative \`a l'infini est nulle, ce qui donne bien les isomorphismes souhait\'es pour $k<(n-3)/2$.

Maintenant, on traite par exemple le cas o\`u $n=2m$ est pair, et o\`u $k=(n-2)/2=m-1$ (le degr\'e $k=(n-3)/2$ lorsque $n$ est impair se traite de fa\c con similaire). On pourra supposer, pour simplifier, qu'on a un seul bout. D'apr\`es ce qui pr\'ec\`ede, on a une injection $$0\ra H^{m-1}_2(M) \buildrel r \over\ra H^{m-1}(D).$$ On va montrer que cette application est surjective. Soit donc $[\alpha ]$ une classe de $H^{m-1}(D)$. Pr\`es de $\partial D$, on \'ecrit $\alpha = \beta (r) +dr \wedge \gamma (r)$, o\`u $\beta (r)$ et $\gamma (r)$ sont consid\'er\'ees comme des formes sur $\partial D$, d\'ependant du param\`etre $r$.  On pose $$\phi (r)=\int _0^r \gamma (t) dt,$$ et on voit facilement que $\alpha -d\phi$ ne contient pas de composante suivant $dr$. De plus, si $d^T$ d\'esigne la diff\'erentielle sur $\partial D$, on a
\begin{eqnarray*}
\frac{\partial }{\partial r}(\alpha -d\phi ) &=& \frac{\partial \beta }{\partial r}-d^T\gamma \\
 &=& d\alpha =0.
\end{eqnarray*}
Donc quitte \`a remplacer $\alpha$ par $\alpha -d\phi$, on pourra supposer que pr\`es de $\partial D$, $\alpha $ est ind\'ependante de $r$, et sans composante normale. Comme l'ext\'erieur de $D$ est diff\'eomorphe \`a un cylindre $]0,\infty [\times \partial D$, on peut prolonger $\alpha$ en dehors de $D$ par $\alpha (r)=\alpha (0)$ pour $r\geq 0$. Ce prolongement, not\'e $\tilde{\alpha }$, reste ferm\'e. De plus, on voit facilement qu'il est de carr\'e int\'egrable, donc $r[\tilde {\alpha }]=[\alpha ]$, et $r$ est bien surjective. 

Pour le degr\'e $n/2$, on utilise l'invariance conforme de la $L^2-$cohomologie en degr\'e moiti\'e : dans ce cas, chaque bout est conform\'ement \'equivalent \`a un cylindre, et la $L^2-$cohomologie d'un cylindre est nulle. On finit ensuite comme dans le point 2. de la preuve du th\'eor\`eme \ref{cohomologie}.\cqfd

\Rq Si dans la formule d'int\'egration par parties, on choisit une fonction de la forme $f(r)$, avec $f'=F$, alors on a 
$$2Hf(\alpha ,\alpha )+\Delta f\vert \alpha \vert ^{2}=F'(2k-n)\vert \alpha \vert ^{2}.$$ 
Alors pour toute forme \`a support compact, le corollaire \ref{IPP1} donne : 
$$\vert \int (2k-n)F'\vert \alpha \vert ^{2}\vert\leq 2 \sup _{support (\alpha )}(F(r)) \left[ \Vert d\alpha \Vert _{L^2}+\Vert \delta \alpha \Vert _{L^2}\right]\Vert \alpha \Vert _{L^2}.$$
Si on suppose que $F'$ est de signe constant (sans s'annuler) et que $F$ est born\'ee, ceci montre que pour $k\neq n/2$, l'op\'erateur $d+\delta$ est non parabolique \`a l'infini en restriction aux $k-$formes, et donc que les espaces de formes harmoniques qui correspondent sont de dimension finie (voir \cite{C3}). En toute rigueur, le r\'esultat de \cite{C3} ne s'applique que si l'op\'erateur $d+\delta$ est non parabolique \`a l'infini globalement, c'est-\`a-dire sur tous les degr\'es. Cependant, la m\^eme preuve s'adapte \`a notre cas. Remarquons aussi que dans \cite{C2}, G. Carron obtient un meilleur r\'esultat si $f$ est concave ($F'<0$) : dans ce cas, et sans autre hypoth\`ese sur $F$, l'op\'erateur $d+\delta $ est non-parabolique \`a l'infini (voir la remarque qui suit \cite[proposition 5.1]{C2}).


 \subsection{Application aux vari\'et\'es conform\'ement compactes} 

Dans cette partie, nous consid\'erons une vari\'et\'e riemannienne $\overline{M}^{n}$ compacte \`a bord, de dimension $n$, munie d'une m\'etrique riemannienne $\overline{g}$ qui est lisse jusqu'au bord. Soit $y : \overline{M}\ra \mathbf{R}_{+}$ une fonction lisse et positive d\'efinissant le bord : $\partial M=y^{-1}(0)$, et $dy\neq 0$ le long de $\partial M$. On munit $M$ de la m\'etrique $g=\overline{g}/y^{2}$ ; cette m\'etrique est compl\`ete et on dit que $\overline{M}$ est une vari\'et\'e conform\'ement compacte. Si de plus la norme de $dy$ par rapport \`a $\overline g$ est constante le long du bord, on dira que $M$ est asymptotiquement hyperbolique. L'exemple typique d'une telle situation est donn\'ee par le mod\`ele du disque hyperbolique, o\`u $\overline{g}$ est la m\'etrique euclidienne, et $y(x)=(1-\vert x\vert ^{2})/2$. Dans \cite{M}, Mazzeo relie la $L^2-$cohomologie de $M$ en degr\'es $k<(n-1)/2$ (respectivement $k>(n+1)/2$) \`a la cohomologie relative (respectivement absolue) de $\overline{M}$. De plus, il d\'etermine compl\`etement le spectre essentiel du laplacien pour tous les degr\'es. Pour ceci, l'auteur utilise des outils de calcul pseudo-diff\'erentiel adapt\'es \`a la g\'eom\'etrie consid\'er\'ee. Ici, nous allons voir comment le th\'eor\`eme \ref{produit tordu} permet de retrouver plus simplement le r\'esultat de Mazzeo concernant la $L^2-$cohomologie :
\begin{cor}[Mazzeo]\label{Mazzeo}
Soit $(M^{n},g)$ une vari\'et\'e conform\'ement compacte. Alors on a les isomorphismes entre espaces vectoriels de dimension finie :
$$H^k_2(M)\simeq \left\{ \begin{array}{ll}
     H^k_c(M), & \textrm{si $k\leq (n-1)/2$,}\\
     H^k(M), & \textrm{si $k\geq(n+1)/2.$}
     \end{array} \right. $$
\end{cor}
\dem On va montrer (et c'est certainement bien connu) qu'\`a l'infini, $M$ est quasi-isom\'etrique \`a un produit tordu $[ 0,\infty [ \times \partial M$ portant la m\'etrique $dr^2+e^{2r}d\theta ^2$ (ici, je remercie encore une fois G. Carron de m'avoir expliqu\'e ce fait). Comme la $L^2-$cohomologie est invariante par quasi-isom\'etries, on aura le r\'esultat gr\^ace au th\'eor\`eme \ref{produit tordu}. D'abord, remarquons que si $\tilde{g}$ et $\tilde{y}$ sont respectivement une autre m\'etrique sur $\overline{M}$ et une autre fonction d\'efinissant le bord, alors les m\'etriques $\overline{g}/y^2$ et $\tilde{g}/\tilde{y}^2$, sont quasi-isom\'etriques. En effet, les m\'etriques $\overline{g}$ et $\tilde{g}$ \'etant d\'efinies sur $\overline{M}$ qui est compacte, elles sont comparables, et comme $y$ et $\tilde{y}$ s'annulent exactement sur $\partial M$ \`a l'ordre z\'ero, $y$ et $\tilde{y}$ sont aussi comparables, et $\overline{g}/y^2$ et $\tilde{g}/\tilde{y}^2$ sont bien quasi-isom\'etriques. Ensuite, si $\overline{g}$ est une m\'etrique donn\'ee, on peut toujours choisir $y$ de sorte que la norme de $dy$ par rapport \`a $\overline{g}$ le long du bord $\partial M$ soit constante : on note $a(\theta )$ cette norme ($\theta$ est la variable sur $\partial M$), et on peut remplacer $y$ par $y/a$ sur un voisinage de $\partial M$. On a alors $d(y/a)=dy/a+yd(1/a)=dy/a$ sur $\partial M$. En r\'esum\'e, on peut supposer que notre m\'etrique conform\'ement compacte est asymptotiquement hyperbolique. Mais pour une telle m\'etrique, la quasi-isom\'etrie avec le produit tordu souhait\'e d\'ecoule par exemple des travaux de Graham (voir \cite{Gr}).\cqfd


\end{document}